\documentclass[11pt]{amsart}

\usepackage{subcaption, booktabs}   
\usepackage{amsmath,amssymb,amsthm,mathtools}
\usepackage{geometry}
\usepackage{graphicx}
\usepackage{booktabs}
\usepackage{hyperref}
\hypersetup{colorlinks=true, linkcolor=blue, citecolor=blue, urlcolor=blue}
\usepackage{enumitem}
\usepackage{tcolorbox}
\tcbuselibrary{breakable,skins}

\newtheorem{theorem}{Theorem}[section]
\newtheorem{proposition}[theorem]{Proposition}
\newtheorem{lemma}[theorem]{Lemma}
\newtheorem{corollary}[theorem]{Corollary}
\theoremstyle{definition}
\newtheorem{definition}[theorem]{Definition}
\newtheorem{assumption}[theorem]{Assumption}
\newtheorem{remark}[theorem]{Remark}

\newtcolorbox{keyidea}[1][]{
  colback=yellow!10, colframe=orange!70!black,
  fonttitle=\bfseries\color{orange!70!black},
  title={Key Idea}, breakable, #1 }

\newcommand{\R}{\mathbb{R}}
\newcommand{\E}{\mathbb{E}}

\newcommand{\RL}[1]{\,{}^{\mathrm{RL}}\!I^{#1}}

\newcommand{\Cder}[1]{\,{}^{\mathrm{C}}\!D^{#1}}       

\newcommand{\rCder}[1]{\,{}^{\mathrm{C}}\!\bar D^{#1}}

\begin{document}

 \title[Time-Fractional Financial  Speculative Bubbles ]{A Time-Fractional Free Boundary Problem for\\ Financial Speculative Bubbles}

\author{Avetik Arakelyan}
\address[A.~Arakelyan]{Institute of Mathematics, SNCO\\Yerevan State University\\ DiGREA lab, CSIE, Yerevan, Armenia}
\email{avetik.arakelyan@ysu.am}

\author{Rafayel Barkhudaryan}
\address[R.~Barkhudaryan]{Institute of Mathematics, SNCO\\ Yerevan State University, Yerevan, Armenia}
\email{barkhudaryan@ysu.am}

\author{Layan El Hajj}
\address[L.~El~Hajj]{Carnegie Mellon University in Qatar, Doha, Qatar}
\email{lelhajj@andrew.cmu.edu}

\author{Henrik Shahgholian}
\address[H.~Shahgholian]{KTH Royal Institute of Technology, Stockholm, Sweden}
\email{henriksh@kth.se}

\begin{abstract}
In this work, we study a time-fractional extension of the Berestycki--Monneau--Scheinkman (BMS) model \cite{BMS14a} for speculative financial bubbles, motivated by the long-range persistence of investor disagreement. Formulating the resale option backward in calendar time up to maturity, we model the speculative premium via a reflected obstacle problem governed by a right Djrbashian--Caputo fractional time derivative of order $\beta\in(0,1]$. The choice $\beta=2-2H$ is motivated by the power-law decay of correlations between increments of fractional Brownian motion with Hurst index $H\in[1/2,1)$, although our analysis treats $\beta$ as an independent parameter.

Our main results establish existence, uniqueness, and regularity of solutions to this non-local-in-time free boundary problem on an unbounded spatial domain. In particular, we prove that the unique viscosity solution is \emph{globally Lipschitz continuous in space and locally $\beta/2$-Hölder continuous in time.} To obtain these results, we develop a viscosity solution framework, establishing equivalent definitions of viscosity solutions, stability under half-relaxed limits, and a comparison principle for semicontinuous sub- and supersolutions. Existence and uniqueness in the uniform linear-growth class follow from Perron's method and explicit sub- and supersolution barriers.

\end{abstract}

\maketitle

\vspace{0.8em}

\vspace{0.4em}

\tableofcontents

\section{Introduction}\label{sec:intro}

Throughout financial history, speculative bubbles have repeatedly accompanied periods of financial instability. While economic studies offer deep insight into the behavioral roots of these phenomena \cite{KA05,Shi15}, capturing their generative mechanics within a formal mathematical model remains an open problem.

The semimartingale framework developed by Jarrow, Protter and their
co-authors \cite{JPS07,JPS10} relates bubbles to strict local martingale
properties of discounted prices under a pricing measure. This approach does
not require pathwise price explosion. To explain the \emph{formation} of bubbles we align
with the equilibrium approach rooted in investor disagreement and limits to
arbitrage. Introduced by Miller \cite{Mil77} and formalized by Harrison and
Kreps \cite{HK78}, this theory posits that bubbles arise from the ``resale
option'': when short-selling is constrained, an asset owner pays a speculative
premium above the fundamental value because they hold an American-style option
to resell to a more optimistic investor. Scheinkman and Xiong \cite{SX03}
transitioned this mechanism into continuous-time diffusions, Chen and Kohn
\cite{CK11} refined it, and Berestycki, Monneau and Scheinkman (BMS)
\cite{BMS14a} reduced the microstructural trading mechanism to a
macroscopic parabolic obstacle PDE whose free boundary represents the optimal
trading threshold between heterogeneous agents.
Related non-local  models with L\'evy jump processes have been studied in \cite{Arakelyan2026arxiv,Arakelyan_2026}.
For the numerical treatment of such non-local free boundary problems we refer to \cite{ABSS19,BJS16}.

The classical BMS framework assumes that the underlying disagreement evolves
along a \emph{memoryless} Markov diffusion. For an Ornstein--Uhlenbeck
specification, stationary autocorrelations decay exponentially. Persistent
sentiment motivates considering models with slower, power-law relaxation.
Fractional Brownian motion (fBm) of Hurst index $H>\tfrac12$ supplies one
mathematical model of long-range dependence \cite{Man68,BHOZ08}. We use this
as modeling motivation, without claiming that every Markov diffusion has
exponentially decaying correlations or that these references establish an
empirical law for sentiment.

In this paper we therefore consider a time-fractional extension of the BMS
obstacle problem. The power-law autocovariance of fBm increments and the
kernel of a \emph{fractional time derivative} suggest the exponent
calibration $\beta=2-2H$. This calibration maps $H\in(\tfrac12,1)$ to
$\beta\in(0,1)$ and $H=\tfrac12$ to $\beta=1$, but does not identify the
underlying stochastic models. Gaussian memory models \cite{KX04} and
inverse-stable time-changed (CTRW) models \cite{MS04,HKU11,GM97,Mag09,SGM00}
must be distinguished. We postulate the fractional obstacle PDE and treat
$\beta$ as an independent order parameter.

We work in \emph{calendar time} $t\in[0,T]$ throughout, with $t=0$ the initial
trading date, $t=T$ maturity, and the terminal condition $u(g,T)=0$ (the resale
option is worthless at maturity). Because the pricing problem is posed backward
from maturity, the correct fractional operator is the \emph{right} (backward)
fractional derivative $\rCder{\beta}_{t}$ on $[t,T]$, whose kernel
$(s-t)^{-\beta}/\Gamma(1-\beta)$ weights the value process between the current
date and maturity, and which at $\beta=1$ reduces to $-\partial_t$, recovering
the classical Brownian BMS backward equation. We study the speculative premium through a \emph{time-fractional} obstacle
problem governed by $\rCder{\beta}_{t}$. The order $\beta$ interpolates between
the classical equation and a family of nonlocal-in-time models.

The presence of the fractional derivative requires a careful adaptation of classical parabolic methods. The operator is \emph{non-local in time}, so the classical
It\^o / Feynman--Kac correspondence and the standard parabolic comparison
proof cannot be transferred without further justification. We work in the viscosity framework of
Crandall--Ishii--Lions \cite{CIL92}, adapted to time-fractional operators
following Giga--Namba \cite{GN17}, Topp--Yangari \cite{TY17} and Namba
\cite{Namba18}.

The problem also connects naturally with   optimal switching and real-options valuation \cite{BayraktarEgami2010,CarmonaLudkovski2008}.
By introducing two reflected regimes and centering the payoff, the scalar obstacle equation can be rewritten as a symmetric system of coupled obstacle equations with constant switching cost $c$; see Section \ref{sec:switching}.

For the particular case,  $\beta=1$, this is the dynamic-programming structure of a classical  two-regime switching problem, whose free boundaries separate continuation from switching.
For $0<\beta<1$, our equation provides a nonlocal-in-time extension of this structure, related to the broader study of fractional Bellman equations \cite{VeretennikovaKolokoltsov2017}. The connection established here is at the level of the governing  equations; a stochastic representation of the fractional switching problem is not assumed.

The rest of the paper is organized as follows.
Section~\ref{sec:model} sets up the model, defines fBm and the fractional
derivative, and formulates the calendar-time time-fractional obstacle PDE.
Section~\ref{sec:visc} defines viscosity solutions, proves the equivalence
of split and global tests, and establishes half-relaxed stability.
Section~\ref{sec:comp} proves the semicontinuous comparison principle.
Section~\ref{sec:exist} constructs the barriers and applies Perron's
method to establish existence.
Section~\ref{sec:regularity} establishes spatial Lipschitz and temporal
H\"older continuity of a viscosity solution. Finally, Section~\ref{sec:switching} describes the connection
with a two-regime optimal switching system and explains the role of
reflection symmetry.

\section{Background and Mathematical Model}\label{sec:model}

We recall the optimal-stopping framework for financial derivative pricing
and the heterogeneous-beliefs model as motivation for our model. Let $\Lambda(g,t)$ be a value function,
$t\in[0,T]$ denoting \emph{calendar time} ($t=0$ the initial date, $t=T$
maturity), driven by a state variable $g_t$. Under risk-neutral valuation with
constant discount rate $r>0$, the general optimal-stopping problem reads
\begin{equation}\label{eqn: opt-stopping}
  \Lambda(g,t)=\sup_{\tau\in[t,T]}
  \E\!\left[e^{-r(\tau-t)}\Phi(g_\tau,\tau)\,\middle|\, g_t=g\right],
\end{equation}
with $\Phi$ the immediate payoff realized upon execution at the stopping time
$\tau$. In the classical Markov setting, the supremum is over admissible
stopping times for the state filtration. For fBm under its natural history
filtration, the current scalar state does not generally suffice for a Markov
dynamic programming principle. Equation~\eqref{eqn: opt-stopping} is not
used below to assert a stochastic derivation of the fractional PDE.

\subsection{Fractional Brownian motion}\label{ssec:fbm}

A possible model for persistent disagreement replaces the Brownian motion
that drives disagreement with a fractional Brownian motion.

\begin{definition}[Fractional Brownian motion]\label{def:fbm}
A \emph{fractional Brownian motion} (fBm) $B^H=\{B^H_t\}_{t\ge0}$ with Hurst
parameter $H\in(0,1)$ is the centred Gaussian process with $B^H_0=0$ a.s.\ and
covariance
\begin{equation}\label{eq:fbm-cov}
  \E\bigl[B^H_t B^H_s\bigr]
  =\tfrac12\bigl(t^{2H}+s^{2H}-|t-s|^{2H}\bigr),\quad t,s\ge0.
\end{equation}
\end{definition}

From \eqref{eq:fbm-cov} we know $B^H$ is \emph{self-similar},
$B^H_{at}\stackrel{d}{=}a^{H}B^H_t$, with $\E[(B^H_t)^2]=t^{2H}$; it has
\emph{stationary increments}; and $B^{1/2}=W$ is standard Brownian motion. The
increment sequence (fractional Gaussian noise) $X_m:=B^H_{m+1}-B^H_m$ has
autocovariance
\begin{equation}\label{eq:fgn-acf}
  \rho_H(m)=\tfrac12\bigl(|m+1|^{2H}-2|m|^{2H}+|m-1|^{2H}\bigr)
  \;\sim\; H(2H-1)\,m^{2H-2}\quad(m\to\infty),
\end{equation}
so for $H>\tfrac12$ the correlations decay proportionally to $m^{-(2-2H)}$ and are not summable. For $H\neq\tfrac12$ the process
$B^H$ is neither Markov nor a semimartingale.

Fractional Brownian motion is closely connected to fractional calculus.
The representation of Mandelbrot and Van Ness integrates from $-\infty$.
A distinct process is Riemann--Liouville fractional Brownian motion
(RL-fBm), initialized at zero \cite{Man68,BHOZ08}:
\begin{equation}\label{eq:mvn}
  Z^H_t=\frac{1}{\Gamma(H+\tfrac12)}\int_0^t (t-s)^{H-1/2}\,dW_s
  =\bigl(\RL{\,H+1/2}\dot W\bigr)(t).
\end{equation}
Its variance is $t^{2H}/[2H\Gamma(H+\tfrac12)^2]$, and its increments are
not generally stationary. Thus $Z^H$ must not be identified with the process
$B^H$ in Definition~\ref{def:fbm}; initiating trading at $t=0$ does not
make these processes equivalent.

A candidate disagreement process is the fractional Ornstein--Uhlenbeck
dynamics driven by fBm, running forward in calendar time,
\begin{equation}\label{eq:fou}
  dg_t=-\rho\,g_t\,dt+\sigma\,dB^H_t,\quad g_0=x.\footnotemark
\end{equation}
\footnotetext{
In classical models driven by standard Brownian motion, It\^o's Lemma utilizes the non-zero quadratic variation of the path to introduce the spatial term $\frac{1}{2}\sigma^2 u_{gg}$. However, for $H > \frac{1}{2}$, fractional Brownian motion has \emph{zero} quadratic variation, hence classical It\^o calculus is inapplicable. Instead, the spatial diffusion term $\frac{1}{2}\sigma^2 u_{gg}$ coupled with the time-fractional derivative arises  as the macroscopic limit of a Continuous Time Random Walk (CTRW) \cite{MS04,Mag09}. In this framework, the spatial price jumps remain standard Gaussian (yielding the local spatial derivatives), while the waiting times between trades follow a  Mittag-Leffler distribution, yielding the non-local fractional time derivative.
}

\subsection{Parameters and the payoff function}\label{ssec:struct}

The intrinsic payoff from reselling the asset depends on the relative
disagreement state $g$ and a deterministic time-decay profile $\alpha$. We
impose $\alpha\in C^1([0,T])$, $\alpha\ge0$, $\alpha(T)=0$, and (reflecting the
dissipation of the resale opportunity towards maturity) $\alpha'(t)\le0$, i.e.\
$\alpha$ is \emph{non-increasing} in calendar time. Such a profile is considered, for instance, in \cite{BMS14a} and is defined as follows
\begin{equation}\label{eq:alpha}
  \alpha(t)=\frac{1-e^{-(r+\lambda)(T-t)}}{r+\lambda},\quad
  \alpha(T)=0,\quad \alpha'(t)=-e^{-(r+\lambda)(T-t)}<0,
\end{equation}
where $\lambda>0$ is the resale (matching) intensity in the classical
interpretation. We take this profile as a prescribed input to the
fractional model; a fractional equilibrium derivation of it is not assumed. Let $c>0$ be the strictly
positive transaction cost. The net payoff realized upon exercise is

\begin{equation}\label{eq:psi}
  \psi(g,t):=g\,\alpha(t)-c.
\end{equation}
This provides the fundamental complementarity condition
\begin{equation}\label{eq:complement}
  \psi(g,t)+\psi(-g,t)=-2c<0,
\end{equation}
For a solution satisfying both reflected obstacle inequalities, their gaps
sum to $2c$, so they cannot both vanish. This secures the disjointness of
the exercise regions.
At maturity, $\psi(g,T)=-c<0$.

\subsection{Fractional derivatives}\label{ssec:fractional}

Long memory is the language of fractional calculus. Because our problem is
posed backward from maturity, we shall need the \emph{right} (backward) fractional
derivative; we recall both one-sided operators.

\begin{definition}[Riemann--Liouville integral and fractional derivatives]
\label{def:fractional}
For $\alpha>0$ the (left) Riemann--Liouville fractional integral of a locally
integrable $f:[0,T]\to\R$ is
$\RL{\alpha}f(t):=\frac{1}{\Gamma(\alpha)}\int_0^t (t-s)^{\alpha-1}f(s)\,ds$.
For $\beta\in(0,1)$ and $f\in C^1([0,T])$, the \emph{left} and \emph{right}
fractional  derivatives are 
\begin{equation}\label{eq:fractional-def}
  \Cder{\beta}_{t}f(t):=\frac{1}{\Gamma(1-\beta)}\int_0^t
  (t-s)^{-\beta}f'(s)\,ds,\quad
  \rCder{\beta}_{t}f(t):=\frac{-1}{\Gamma(1-\beta)}\int_t^T
  (s-t)^{-\beta}f'(s)\,ds.
\end{equation}
For absolutely continuous $f$, the integral definitions are used at points
where the weighted integrals are finite. By convention
$\rCder{1}=-\tfrac{d}{dt}$. No order-zero convention is needed; the
order-zero fractional limit is $f(t)-f(T)$, not generally $f(t)$. The right
derivative \eqref{eq:fractional-def} integrates over the \emph{future} $[t,T]$,
which is the natural memory for a value posed backward from maturity. It admits
the equivalent regularized (Marchaud-type) form \cite{kilbas2006theory}, valid for smooth $f$ and, as a difference operator, whenever the integral
is absolutely convergent,

\begin{equation}\label{eq:marchaud}
  \rCder{\beta}_{t}f(t)
  =\frac{f(t)-f(T)}{\Gamma(1-\beta)\,(T-t)^{\beta}}
  +\frac{\beta}{\Gamma(1-\beta)}\int_t^T
  \frac{f(t)-f(s)}{(s-t)^{1+\beta}}\,ds .
\end{equation}
Under the time reversal $\tau=T-t$, the right derivative in $t$ equals the left
derivative in $\tau$, i.e.  $\rCder{\beta}_{t}f(t)=\Cder{\beta}_{\tau}\tilde f(\tau)$
with $\tilde f(\tau)=f(T-\tau)$.
\end{definition}

\begin{remark}[Historical note on \eqref{eq:fractional-def}]\label{rem:history}
The fractional operator in which the ordinary derivative is taken \emph{before}
the fractional integration is widely attributed to M.~Caputo, who introduced
it in $1967$ in the study of dissipative media \cite{Caputo67}. It was, however,
independently and essentially simultaneously introduced by the Armenian
mathematicians M.~M.~Djrbashian (Dzhrbashian) and A.~B.~Nersessian in their
$1968$ paper on the Cauchy problem for fractional differential equations
\cite{DN68}, building on Djrbashian's earlier work on fractional
integro-differentiation and Mittag-Leffler-type entire functions
\cite{Djrbashian66}. For this reason the operator is increasingly called the
\emph{Djrbashian--Caputo} fractional derivative \cite{Rogosin18}. The feature
decisive below that it admits classical (integer-order) terminal data, unlike
the Riemann--Liouville derivative was already central to the
Djrbashian--Nersessian treatment.
\end{remark}

\subsection{Fractional obstacle PDE with fractional time derivative}
\label{ssec:derivation}

We formulate a  time-fractional governing equation motivated by
power-law memory. The spatial part, describing mean reversion and diffusion, is
taken from the classical model. This is a modeling postulate, not a
consequence of the fBm stochastic differential equation \eqref{eq:fou}.

\paragraph{The memory order.}
The fractional Gaussian noise driving \eqref{eq:fou} has the autocovariance
\eqref{eq:fgn-acf}, decaying as $\rho_H(m)\sim m^{-(2-2H)}$ for
$H>\tfrac12$. Fractional relaxation equations have Mittag-Leffler
relaxation with power-law tails \cite{Mainardi96}. The relation
$\beta=2-2H$ can therefore be considered as an exponent calibration, but
we leave $\beta\in(0,1]$ independent in the analysis.

The distinction is already visible without mean reversion. The variance of
$\sigma B^H_t$ is $\sigma^2t^{2H}$, whereas the fractional heat equation
$\Cder{\beta}_t p=(\sigma^2/2)p_{gg}$ with initial unit mass at zero
has second moment $\sigma^2t^\beta/\Gamma(1+\beta)$. These powers differ
when $H>\tfrac12$ and $\beta=2-2H$. Inverse-stable time changes and fBm
thus require distinct evolution equations \cite{MS04,HKU11}.

\paragraph{The obstacle problem.}
Upon exercise the seller receives the immediate cash flow $\psi(g,t)$ plus the
continuing value from the buyers' inverse perspective $-g$, so the payoff is
coupled to the solution itself.

In a classical Markov interpretation, applying this to
\eqref{eqn: opt-stopping} motivates the bubble valuation equation
\[
\Lambda_{\text{Bubble}}(g, t) = \sup_{\tau \in [t, T]} \mathbb{E}\!\left[ e^{-r(\tau-t)} \left[ \Lambda_{\text{Bubble}}(-g_\tau, \tau) + \psi(g_\tau, \tau) \right]^+ \;\Big|\; g_t = g \right].
\]
\begin{definition}[Time-fractional bubble obstacle PDE in calendar time]
\label{thm:pde}
For $\beta\in(0,1]$, we assume that the bubble premium $u(g,t)$ is
described by the non-local-in-time obstacle problem
\begin{equation}\label{eq:pde}
  \min\Bigl\{\,F[u](g,t),\;\; u(g,t)-u(-g,t)-\psi(g,t)\Bigr\}=0,
  \quad u(g,T)=0,
\end{equation}
where, with the \emph{right} (backward) fractional derivative in calendar time
$t\in[0,T]$,
\begin{equation}\label{eq:operator}
  F[u]:=\rCder{\beta}_{t}u+\rho g\,u_g-\tfrac{\sigma^2}{2}u_{gg}+r u,
  \quad \beta\in(0,1].
\end{equation}
For $\beta=1$ ($H=\tfrac12$), $\rCder{1}_t=-\partial_t$ and \eqref{eq:pde}
becomes the classical bubble obstacle problem \cite{BMS14a}.
\end{definition}

For classical, memoryless disagreement ($H=\tfrac12$, $\beta=1$), the
continuation value solves the backward parabolic equation
$-\partial_t u=\mathcal A u-ru$, where
$\mathcal A u=-\rho g u_g+\tfrac{\sigma^2}{2}u_{gg}$. Moving all terms
to the left yields the standard BMS operator. In the fractional model,
$-\partial_t$ is replaced by $\rCder{\beta}_t$ as a prescribed extension.
The obstacle problem partitions the domain into two regimes:
\begin{enumerate}
\item \emph{The continuation region:} the value is strictly greater than the
immediate execution payoff, $u(g,t)>u(-g,t)+\psi(g,t)$, and $F[u]=0$.
\item \emph{The exercise region:} the value equals the obstacle,
$u(g,t)=u(-g,t)+\psi(g,t)$, and the model requires $F[u]\ge0$.
Writing $B[u](g,t):=u(g,t)-u(-g,t)-\psi(g,t)$, we have
$B[u](g,t)+B[u](-g,t)=2c$. Contact at $g$ therefore gives
$B[u](-g,t)=2c>0$. This does not require $\psi(g,t)\ge0$.
\end{enumerate}
The inequalities hold complementarily at each point, with equality in at
least one branch. The differential inequality on the exercise region is
part of the postulated problem: a pointwise contact identity cannot be
differentiated with a time-nonlocal operator without information along the
whole interval $[t,T]$.

At maturity $t=T$, the resale option expires. Since $\alpha(T)=0$, the
immediate exercise payoff is $\psi(g,T)=-c$. The terminal condition is
therefore $u(g,T)=0$. We seek solutions in the uniform linear-growth class
$|u(g,t)|\le C(1+|g|)$ on $\R\times[0,T]$.

A stochastic representation would require an additional choice of clock,
filtration, and admissible stopping strategies. In particular,
$\Cder{\beta}_s q+rq=0$, $q(0)=1$, has solution
$q(s)=E_\beta(-rs^\beta)$, with
$E_\beta(z)=\sum_{k\ge0}z^k/\Gamma(1+\beta k)$; deterministic
calendar-time discounting instead gives $e^{-rs}$. Therefore the
fractional equation is not identified here with the preceding stopping
formula for fBm. Coefficients in the postulated equation may be understood
as nondimensionalized; their classical calendar-time interpretations are
not automatically preserved for $\beta<1$.

\begin{assumption}\label{ass:standing}
Throughout the paper we assume $r>0$, $\rho\ge0$, $\sigma>0$, $c>0$, $T>0$,
$\beta\in(0,1]$ (equivalently $H\in[\tfrac12,1)$ under the optional
parametrization $\beta=2-2H$).
\end{assumption}

\section{Viscosity Solutions}\label{sec:visc}

Classical solutions require $u\in C^2$ in space and $\rCder{\beta}_t u$ to exist
pointwise, which is not assumed across the free boundary. We use the viscosity
framework of Crandall--Ishii--Lions \cite{CIL92}, adapted to the fractional operator
following \cite{GN17,TY17,Namba18}. At a normalized global upper contact,
positivity of the Marchaud kernel gives
$\rCder{\beta}_t\varphi(g_0,t_0)\le\rCder{\beta}_t u(g_0,t_0)$
whenever both expressions exist. Global-test formulations with a suitable
test class are legitimate \cite{Namba18}. For the proofs below we explicitly
write the actual solution in the distant-time part of the operator.

\begin{definition}\label{def:testop}
Let $Q:=\R\times[0,T]$ and $Q_0:=\R\times[0,T).$ We set 
\[
\mathcal T:=C(Q)\cap C^{2,1}(Q_0).
\]
For $(g_0,t_0)\in Q_0$ and $\varphi\in\mathcal T$, define
\begin{equation}\label{eq:testop}
  F_\varphi(g_0,t_0;u):=\rCder{\beta}_{t}\varphi(g_0,t_0)
  +\rho g_0\,\varphi_g(g_0,t_0)
  -\tfrac{\sigma^2}{2}\varphi_{gg}(g_0,t_0)
  +r\,u(g_0,t_0).
\end{equation}
For $0<\beta<1$, the operator on tests is defined by the Marchaud
formula \eqref{eq:marchaud}; for $\beta=1$, it is $-\varphi_t$.
Time derivatives at $t=0$ are right-sided. The zero-order term
carries the actual $u$, which agrees with the test value at a
normalized contact. No spatial growth restriction is imposed on tests.
\end{definition}

For $0<\beta<1$ and $0<\delta<T-t_0$, define the split operator
\begin{align*}
D_\delta[u,\varphi](g_0,t_0)
&=\frac{u(g_0,t_0)-u(g_0,T)}{\Gamma(1-\beta)(T-t_0)^\beta}\\
&+\frac{\beta}{\Gamma(1-\beta)}\int_0^\delta
\frac{\varphi(g_0,t_0)-\varphi(g_0,t_0+\tau)}{\tau^{1+\beta}}\,d\tau\\
&+\frac{\beta}{\Gamma(1-\beta)}\int_\delta^{T-t_0}
\frac{u(g_0,t_0)-u(g_0,t_0+\tau)}{\tau^{1+\beta}}\,d\tau,
\end{align*}
and let $F_{\varphi,\delta}(g_0,t_0;u)$ be \eqref{eq:testop} with
$\rCder{\beta}_t\varphi$ replaced by $D_\delta[u,\varphi]$.
For $\beta=1$ set $F_{\varphi,\delta}=F_\varphi$, using $-\varphi_t$.

\begin{definition}[Viscosity solutions of \eqref{eq:pde}]\label{def:visc}
All candidates are locally bounded on the closed strip. Contacts are \textit{normalized}\footnote{At the contact point,
$\varphi(g_0,t_0)=u(g_0,t_0)$ for subsolutions, and
$\varphi(g_0,t_0)=v(g_0,t_0)$ for supersolutions. Thus the corresponding
local maximum of $u-\varphi$, or local minimum of $v-\varphi$, equals
zero. This can always be achieved by adding a constant to $\varphi$,
which leaves its derivatives and the increments in the split
fractional operator unchanged.}
by equality of the function and test value. The test is $C^{2,1}$
in its contact neighborhood; $\delta$ is admissible when
$0<\delta<T-t_0$ and the segment $\{(g_0,t_0+\tau):0\le \tau\le\delta\}$
lies in that neighborhood. For every $\beta\in(0,1]$, contacts at $t=0$ are taken relative to the closed strip. For $\beta=1$, the temporal term at $t=0$
is the right-sided derivative $-\varphi_t(g_0,0)$.
\leavevmode
\begin{itemize}[leftmargin=1.4em]
\item An upper-semicontinuous (USC) $u:\R\times[0,T]\to\R$ with $u(\cdot,T)\le0$
is a \textbf{\emph{viscosity subsolution}} if for every $\varphi\in C^{2,1}$ and every
normalized \emph{local} maximum point  $(g_0,t_0)$, $t_0<T$, of $u-\varphi$, the following inequality holds for every admissible $\delta$
\begin{equation}\label{eq:sub}
  \min\Bigl\{F_{\varphi,\delta}(g_0,t_0;u),\;
  u(g_0,t_0)-u(-g_0,t_0)-\psi(g_0,t_0)\Bigr\}\le0.
\end{equation}
\item A lower-semicontinuous (LSC) $v:\R\times[0,T]\to\R$ with $v(\cdot,T)\ge0$
satisfying $v(g,t)-v(-g,t)-\psi(g,t)\ge0$ for $(g,t)\in Q_0$
is a \textbf{\emph{viscosity supersolution}} if for every $\varphi\in C^{2,1}$ and every
normalized local minimum point $(g_0,t_0)$, $t_0<T$, of $v-\varphi$, the following inequality holds for every admissible $\delta$
\begin{equation}\label{eq:super}
  \min\Bigl\{F_{\varphi,\delta}(g_0,t_0;v),\;
  v(g_0,t_0)-v(-g_0,t_0)-\psi(g_0,t_0)\Bigr\}\ge0.
\end{equation}
\item $u$ is a \textbf{\emph{viscosity solution}} if $u^*$ is a subsolution and $u_*$ a
supersolution.
\end{itemize}
\end{definition}

\begin{lemma}[Equivalence of split and global tests]\label{lem:global-split}
Let $u\in USC(Q)$ be locally bounded and satisfy
$u(g,t)\le C(1+|g|)$ and $u(\cdot,T)\le0$.
Then $u$ is a subsolution in Definition~\ref{def:visc} if and
only if, at every normalized global maximum $(g_0,t_0)$ of
$u-\varphi$, with $\varphi\in\mathcal T$ and $t_0<T$,
\[
\min\{F_\varphi(g_0,t_0;u),B[u](g_0,t_0)\}\le0.
\]
Likewise, let $v\in LSC(Q)$ be locally bounded and satisfy
$v(g,t)\ge-C(1+|g|)$ and $v(\cdot,T)\ge0$.
Then $v$ is a supersolution in Definition~\ref{def:visc} if and
only if, at every normalized global minimum $(g_0,t_0)$ of
$v-\varphi$, with $\varphi\in\mathcal T$ and $t_0<T$,
\[
\min\{F_\varphi(g_0,t_0;v),B[v](g_0,t_0)\}\ge0.
\]
Here $B[w](g,t)=w(g,t)-w(-g,t)-\psi(g,t)$, and contacts at
$t=0$ use the convention in Definition~\ref{def:visc}.
\end{lemma}

\begin{proof}
We first take $0<\beta<1$ and write $L:=T-t_0$.

Suppose $u$ satisfies the split definition and $\varphi$ is a
normalized global upper test at $(g_0,t_0)$. Then $\varphi\ge u$
on $Q$, and, for every sufficiently small admissible $\delta$,
\begin{align*}
&F_{\varphi,\delta}(g_0,t_0;u)-F_\varphi(g_0,t_0;u)\\
&=\frac{\varphi(g_0,T)-u(g_0,T)}{\Gamma(1-\beta)L^\beta}
 +\frac{\beta}{\Gamma(1-\beta)}\int_\delta^L
 \frac{\varphi(g_0,t_0+\tau)-u(g_0,t_0+\tau)}{\tau^{1+\beta}}\,d\tau
 \ge0.
\end{align*}
Therefore,
\[
\min\{F_\varphi,B[u]\}
\le\min\{F_{\varphi,\delta},B[u]\}\le0
\quad\text{at }(g_0,t_0).
\]
At a normalized global lower contact with $v$, the same calculation
reverses the operator inequality and gives
\[
\min\{F_\varphi,B[v]\}
\ge\min\{F_{\varphi,\delta},B[v]\}\ge0.
\]

Conversely, suppose that $u$ satisfies the global-test inequalities.
Let $\varphi$ be a normalized local upper test at $(g_0,t_0)$
and fix an admissible $\delta$. Choose a contact neighborhood $U$
such that $\varphi\ge u$ on $U$ and
\[
\{(g_0,t_0+\tau):0\le \tau\le\delta\}\subset U.
\]

We require smooth majorants $U_n\ge u$ with $U_n\to u$
pointwise on $Q$, locally uniformly bounded in $n$. To construct
them, for integers $n>C$ set, on $\R^2$,
\[
H_n(x):=\sup_{y\in Q}\{u(y)-n|x-y|\},
\quad
U_n(x):=\int_{\vert{}z\vert{} \le n^{-2}} H_n(x - z)\zeta_{n^{-2}}(z)\,dz+\frac1n,
\]
where $\zeta_\varepsilon$ is a nonnegative smooth mollifier
supported in the ball of radius $\varepsilon$ and of integral one.
The growth bound makes $H_n$ finite and $n$-Lipschitz.
Upper semicontinuity gives $H_n\to u$ on $Q$, and
\[
u\le H_n\le U_n\le H_n+\frac2n.
\]
These inequalities also give the required local uniform bounds.

Choose $\varepsilon_n\to0$ with
$\delta+\varepsilon_n<L$, and smooth cutoffs on $Q$,
$0\le\chi_n\le1$, supported in a fixed compact subset of $U$
away from $t=T$, such that
$\chi_n=1$ near $(g_0,t_0)$ and
\[
\chi_n(g_0,t_0+\tau)=
\begin{cases}
1,&0\le \tau\le\delta,\\
0,&\delta+\varepsilon_n\le \tau\le L.
\end{cases}
\]
Such cutoffs exist because the closed near-time segment is
contained in the contact neighborhood. Define
\[
\varphi_n:=\chi_n\varphi+(1-\chi_n)U_n,
\]
extending $\chi_n\varphi$ by zero outside its support.
Then $\varphi_n\in\mathcal T$ and
\[
\varphi_n\ge u\ \text{on }Q,\quad
\varphi_n(g_0,t_0)=u(g_0,t_0),\quad
\varphi_n=\varphi\ \text{near }(g_0,t_0).
\]
Moreover,
\begin{align*}
\varphi_n(g_0,t_0+\tau)&=\varphi(g_0,t_0+\tau),
\quad 0\le \tau\le\delta,\\
\varphi_n(g_0,t_0+\tau)&\to u(g_0,t_0+\tau),
\quad\delta<\tau\le L.
\end{align*}
The near-time integral is therefore unchanged, while dominated
convergence applies to the distant-time integral, since
$\tau^{-1-\beta}$ is integrable on $[\delta,L]$ and the tests are
uniformly bounded along this fixed history. At the terminal point,
$\varphi_n(g_0,T)=U_n(g_0,T)\to u(g_0,T)$. Hence
\begin{align*}
\rCder{\beta}_t\varphi_n(g_0,t_0)
&\to
\frac{u(g_0,t_0)-u(g_0,T)}{\Gamma(1-\beta)L^\beta}\\
&+\frac{\beta}{\Gamma(1-\beta)}\int_0^\delta
\frac{\varphi(g_0,t_0)-\varphi(g_0,t_0+\tau)}{\tau^{1+\beta}}\,d\tau\\
&+\frac{\beta}{\Gamma(1-\beta)}\int_\delta^L
\frac{u(g_0,t_0)-u(g_0,t_0+\tau)}{\tau^{1+\beta}}\,d\tau\\
&=D_\delta[u,\varphi](g_0,t_0).
\end{align*}
The local derivatives agree at the contact point, so
\[
F_{\varphi_n}(g_0,t_0;u)
\longrightarrow F_{\varphi,\delta}(g_0,t_0;u).
\]
Applying the global inequality to $\varphi_n$ and keeping
$B[u](g_0,t_0)$ fixed yields
\begin{align*}
0\ge\lim_{n\to\infty}
\min\{F_{\varphi_n}(g_0,t_0;u),B[u](g_0,t_0)\}
=\min\{F_{\varphi,\delta}(g_0,t_0;u),B[u](g_0,t_0)\}.
\end{align*}
This proves the split subsolution inequality for every admissible
$\delta$.

For a local lower contact with $v$, apply the same approximation
to $-v$ and negate it to obtain smooth minorants $V_n\le v$,
$V_n\to v$. The tests
$\chi_n\varphi+(1-\chi_n)V_n$ touch $v$ globally from below at
the same point. The preceding limit then gives
\[
\min\{F_{\varphi,\delta}(g_0,t_0;v),B[v](g_0,t_0)\}\ge0.
\]

For $\beta=1$, the two temporal operators agree locally.
The same gluing preserves $\varphi_t$, $\varphi_g$ and
$\varphi_{gg}$ at the contact point, including the right-sided
time derivative at $t_0=0$, so it gives the same equivalence of
the contact inequalities.

It remains to recover the pointwise obstacle constraint for a
global-test supersolution. Fix $p_0=(g_0,t_0)$ with $t_0<T$,
and let $p_n=(g_n,t_n)$ minimize
\[
v(p)+n|p-p_0|^2\quad\text{over }p\in Q.
\]

Since $p_n = (g_n, t_n)$ is the global minimizer over $Q$, evaluating the penalized function at $p_n$ gives a value less than or equal to its value at $p = p_0$
\begin{equation*}
    v(p_n) + n\vert p_n - p_0 \vert^2 \le v(p_0) + n\vert p_0 - p_0 \vert^2 = v(p_0).
\end{equation*}
Using the lower growth bound $v(p_n) \ge -C(1+\vert g_n \vert)$, and Young’s inequality gives
\begin{align*}
    n\vert p_n - p_0 \vert^2  &\le v(p_0) - v(p_n) \le v(p_0) + C(1+\vert g_n \vert) \le v(p_0) + C\left(1 + \vert{}g_0\vert{} + |p_n-p_0|\right)\\&
    \le v(p_0) + C\left(1 + \vert{}g_0\vert{}\right) + \frac{n}{2}\vert p_n - p_0 \vert^2 +\frac{C^2}{2n}.
\end{align*}
 This leads to 
 \[
 n\vert p_n - p_0 \vert^2 \le 2(v(p_0) + C\left(1 + \vert{}g_0\vert{}\right))+\frac{C^2}{n}.
 \]
The last inequality implies that $n\vert p_n - p_0 \vert^2$ cannot blow up as $n \to \infty$. This forces $\vert p_n - p_0 \vert \to 0$, and therefore $p_n \to p_0$ as $n \to \infty$.
   Furthermore, since $n\vert p_n - p_0 \vert^2 \ge 0$, we get $v(p_n) \le v(p_0)$ for all $n$, hence
    \begin{equation*}
        \limsup_{n \to \infty} v(p_n) \le v(p_0).
    \end{equation*}
    Conversely, because $v$ is lower semicontinuous at $p_0$ and $p_n \to p_0$, we have
    \begin{equation*}
        \liminf_{n \to \infty} v(p_n) \ge v(p_0).
    \end{equation*}
    Combining the two inequalities forces exact convergence 
\[
p_n\longrightarrow p_0,\quad
v(p_n)\longrightarrow v(p_0).
\]
The smooth function
\[
q_n(p):=v(p_n)+n|p_n-p_0|^2-n|p-p_0|^2
\]
touches $v$ globally from below at $p_n$. For large $n$ we have
$t_n<T$, so the global supersolution inequality gives
$B[v](g_n,t_n)\ge0$. Lower semicontinuity at the reflected point
and continuity of $\psi$ imply
\begin{align*}
B[v](g_0,t_0)
\ge\limsup_{n\to\infty}
B[v](g_n,t_n)
\ge0.
\end{align*}
Thus the pointwise obstacle constraint also holds on
$\R\times[0,T)$. The terminal inequalities are identical in
both formulations, completing the proof.
\end{proof}

By Lemma~\ref{lem:global-split}, the stability argument and
Perron's construction below may use global tests in $\mathcal T$;
their subsolution and supersolution conclusions are understood
in Definition~\ref{def:visc}.

In viscosity theory, standard pointwise limits need not preserve
semicontinuity. To overcome this difficulty, the stability theory developed
by Ishii \cite{Ishii1987} and Barles and Perthame
\cite{BarlesPerthame1987} uses the half-relaxed limits
($\limsup^*$ and $\liminf_*$). By allowing the spatial-temporal point to
vary together with the family index, these limits preserve the appropriate
semicontinuity. Under the hypotheses below, they also preserve the
subsolution and supersolution inequalities, providing the stability
properties needed in Perron's method.

\begin{proposition}[Stability]\label{prop:stab}
Let $\{u_\alpha\}_{\alpha>0}$ be a family of viscosity subsolutions of
\eqref{eq:pde}, locally uniformly bounded on $\R\times[0,T]$, in the
global-test formulation with tests in
\[
C(Q)\cap C^{2,1}(Q_0).
\]
For $0<\beta<1$, the operator on tests is understood through
\eqref{eq:marchaud}; for $\beta=1$, it is $-\partial_t$, with a
one-sided derivative at $t=0$.
Suppose that there exists a continuous function $\mathsf M $ of at most linear
growth such that
\[
u_\alpha(g,t)\le \mathsf{M} (g,t)
\quad\text{on }Q,\quad
\mathsf M(g,T)\le0.
\]
Then the relaxed upper envelope
\[
\bar u(g,t):=\limsup_{\substack{\alpha\to0\\(y,s)\to(g,t)}}
u_\alpha(y,s),\quad (y,s)\in\R\times[0,T],
\]
is a viscosity subsolution of \eqref{eq:pde} and
$\bar u(\cdot,T)\le0$.

The symmetric statement holds for the lower half-relaxed limit of a
locally uniformly bounded family of viscosity supersolutions, provided
that the family is bounded below by a continuous function $\mathsf m$ of at most
linear growth satisfying $\mathsf m(g,T)\ge0$.
The same conclusions hold for the upper semicontinuous envelope of the
pointwise supremum of subsolutions and the lower semicontinuous envelope
of the pointwise infimum of supersolutions, under the corresponding bounds.
\end{proposition}

\begin{proof}
Let $\varphi\in C(Q)\cap C^{2,1}(Q_0)$
and let $(g_0,t_0)$, with $t_0<T$, be a
strict global maximum of $\bar u-\varphi$. Without loss of generality,
assume
\[
\bar u(g_0,t_0)=\varphi(g_0,t_0),\quad
\bar u-\varphi\le0\quad\text{on }\R\times[0,T].
\]
Since the time interval $[0,T]$ is compact, it remains to control the
maximizing points at spatial infinity. Choose a purely spatial function
$f\in C^\infty(\R)$, $f\ge0$, which vanishes on a neighborhood of $g_0$
and satisfies, for all sufficiently large $|g|$,
\[
f(g)\ge g^2+
\max_{s\in[0,T]}\bigl(\mathsf M(g,s)-\varphi(g,s)\bigr)_+.
\]
Such a function can be obtained by taking a smooth majorant away from
that neighborhood. Set $\widehat\varphi:=\varphi+f$. Then
\[
u_\alpha(g,t)-\widehat\varphi(g,t)
\le \mathsf M(g,t)-\varphi(g,t)-f(g)
\le-g^2
\]
for large $|g|$, uniformly in $\alpha$ and $t$. Thus the maxima are
attained in a bounded spatial region. Moreover,
\[
\rCder{\beta}_t f(g)=0,
\quad
f=f'=f''=0\quad\text{near }g_0,
\]
so the addition of $f$ leaves $F_\varphi$ unchanged near the contact point.
The strict global maximum of $\bar u-\widehat\varphi$ remains
$(g_0,t_0)$, with value zero.

By the definition of the half-relaxed limit, there exist
$\alpha_n\to0$ and $(\hat g_n,\hat t_n)\to(g_0,t_0)$ such that
\[
u_{\alpha_n}(\hat g_n,\hat t_n)\longrightarrow\bar u(g_0,t_0).
\]
Choose $(g_n,t_n)$ to be a global maximum point of
$u_{\alpha_n}-\widehat\varphi$, and set
\[
a_n:=u_{\alpha_n}(g_n,t_n)-\widehat\varphi(g_n,t_n).
\]
Then
\[
a_n\ge u_{\alpha_n}(\hat g_n,\hat t_n)
-\widehat\varphi(\hat g_n,\hat t_n)\longrightarrow0.
\]
The uniform coercive bound above keeps $(g_n,t_n)$ in a fixed compact
set for all sufficiently large $n$. Along any subsequence for which
$(g_n,t_n)\to(\hat g,\hat t)$, the half-relaxed-limit inequality yields
\[
0\le\liminf_{n\to\infty}a_n
\le\limsup_{n\to\infty}a_n
\le\bar u(\hat g,\hat t)-\widehat\varphi(\hat g,\hat t)
\le0.
\]
Strictness of the contact therefore implies
$(\hat g,\hat t)=(g_0,t_0)$. Therefore,
\[
(g_n,t_n)\to(g_0,t_0),\quad a_n\to0,
\]
and
\begin{equation}\label{eq:limit_u}
u_{\alpha_n}(g_n,t_n)
=\widehat\varphi(g_n,t_n)+a_n
\longrightarrow\bar u(g_0,t_0).
\end{equation}

For all sufficiently large $n$, $t_n<T$ and $g_n$ lies in the region
where $f=f'=f''=0$. The normalized upper test for $u_{\alpha_n}$ is
$\widehat\varphi+a_n$. Since the additive constant has zero fractional
derivative and $F_\varphi$ uses the actual candidate in its zero-order term,
\[
F_{\widehat\varphi+a_n}(g_n,t_n;u_{\alpha_n})
=F_\varphi(g_n,t_n;u_{\alpha_n}).
\]
Hence the subsolution inequality \eqref{eq:sub} gives
\begin{equation}\label{eq:pre_limit_sub}
\min\Bigl\{F_\varphi(g_n,t_n;u_{\alpha_n}),\;
 u_{\alpha_n}(g_n,t_n)-u_{\alpha_n}(-g_n,t_n)-\psi(g_n,t_n)
\Bigr\}\le0.
\end{equation}
For every $n$, either the differential term is nonpositive or the obstacle
term is nonpositive. Passing to a subsequence, one of these alternatives
holds for every $n$. We evaluate both cases.

Suppose
$F_\varphi(g_n,t_n;u_{\alpha_n})\le0$ along the subsequence. Expanding,
\[
\rCder{\beta}_t\varphi(g_n,t_n)
+\rho g_n\varphi_g(g_n,t_n)
-\tfrac{\sigma^2}{2}\varphi_{gg}(g_n,t_n)
+r u_{\alpha_n}(g_n,t_n)\le0.
\]
The spatial derivatives are continuous, so
\[
\varphi_g(g_n,t_n)\to\varphi_g(g_0,t_0),\quad
\varphi_{gg}(g_n,t_n)\to\varphi_{gg}(g_0,t_0).
\]
Furthermore, the right fractional operator of the test is continuous
at $(g_0,t_0)$. Indeed, for $0<\beta<1$, let $\mathcal K$ be a compact
neighborhood of $(g_0,t_0)$ relative to $\R\times[0,T)$.
For sufficiently small $\delta>0$, local boundedness of $\varphi_t$ gives
\[
\sup_{(g,t)\in K}\int_0^\delta
\frac{|\varphi(g,t)-\varphi(g,t+\tau)|}{\tau^{1+\beta}}\,d\tau
\le\frac{C_K\delta^{1-\beta}}{1-\beta}
\longrightarrow0\quad(\delta\to0).
\]
For fixed $\delta$, the terminal term and the integral over
$\tau\in[\delta,T-t]$ in \eqref{eq:marchaud} are continuous on $K$,
by continuity of $\varphi$ up to $T$. The uniform estimate above therefore
implies
\[
\rCder{\beta}_t\varphi(g_n,t_n)
\longrightarrow\rCder{\beta}_t\varphi(g_0,t_0).
\]
For $\beta=1$, the same conclusion follows from continuity of
$-\varphi_t$. Combining these limits with \eqref{eq:limit_u}, we obtain
\[
F_\varphi(g_0,t_0;\bar u)\le0.
\]

Suppose instead that
\[
u_{\alpha_n}(g_n,t_n)-u_{\alpha_n}(-g_n,t_n)-\psi(g_n,t_n)\le0
\]
along the subsequence. At the conjugate point, the definition of the
relaxed upper envelope gives
\[
\limsup_{n\to\infty}u_{\alpha_n}(-g_n,t_n)
\le\bar u(-g_0,t_0).
\]
Equivalently,
\[
\liminf_{n\to\infty}\bigl[-u_{\alpha_n}(-g_n,t_n)\bigr]
\ge-\bar u(-g_0,t_0).
\]
Using \eqref{eq:limit_u} and continuity of $\psi$, we obtain
\begin{align*}
0
&\ge\liminf_{n\to\infty}
\Bigl[u_{\alpha_n}(g_n,t_n)-u_{\alpha_n}(-g_n,t_n)-\psi(g_n,t_n)\Bigr]\\
&=\bar u(g_0,t_0)
 +\liminf_{n\to\infty}\bigl[-u_{\alpha_n}(-g_n,t_n)\bigr]
 -\psi(g_0,t_0)\\
&\ge\bar u(g_0,t_0)-\bar u(-g_0,t_0)-\psi(g_0,t_0).
\end{align*}
Thus the obstacle inequality holds at the limit point.

In either case,
\[
\min\Bigl\{F_\varphi(g_0,t_0;\bar u),\;
\bar u(g_0,t_0)-\bar u(-g_0,t_0)-\psi(g_0,t_0)\Bigr\}\le0.
\]
Non-strict global contacts are treated by replacing $\varphi$ with
$\varphi+\varepsilon[(g-g_0)^4+(t-t_0)^4]$ and then letting
$\varepsilon\to0$. The added term and its first two spatial
derivatives vanish at the contact point, and its fractional contribution tends
to zero. Hence $\bar u$ is a viscosity subsolution on $\R\times[0,T)$.

Finally, continuity of the common majorant supplies the terminal condition:
\[
\bar u(g,T)
\le\limsup_{\substack{\alpha\to0\\(y,s)\to(g,T)}}\mathsf M(y,s)
=\mathsf M(g,T)\le0.
\]

For supersolutions, let $\underline v:=\liminf_*v_\alpha$. The same
localization argument, subtracting a suitable purely spatial penalty from
a global lower test, produces contact points with
\[
v_{\alpha_n}(g_n,t_n)\to\underline v(g_0,t_0),\quad
\liminf_{n\to\infty}v_{\alpha_n}(-g_n,t_n)
\ge\underline v(-g_0,t_0).
\]
Both branches of each supersolution inequality are nonnegative. Therefore,
\begin{align*}
F_\varphi(g_0,t_0;\underline v)
&=\lim_{n\to\infty}F_\varphi(g_n,t_n;v_{\alpha_n})\ge0,\\
\underline v(g_0,t_0)-\underline v(-g_0,t_0)-\psi(g_0,t_0)
&\ge\limsup_{n\to\infty}
\Bigl[v_{\alpha_n}(g_n,t_n)-v_{\alpha_n}(-g_n,t_n)-\psi(g_n,t_n)\Bigr]
\ge0.
\end{align*}
The continuous minorant gives
$\underline v(g,T)\ge \mathsf m(g,T)\ge0$.

For the envelopes of pointwise suprema or infima, choose indices that
approximate the corresponding envelope, without requiring
$\alpha_n\to0$. The contact selection and reflected-value estimates
above remain valid, proving the last assertion.
\end{proof}

\section{Comparison Principle}\label{sec:comp}

The comparison argument is formulated directly for semicontinuous
functions, following the regularization strategy for fractional equations in
\cite{Namba18}. Under the reversal $\tau=T-t$, the right fractional
operator on smooth tests becomes the left fractional operator. For the
present affine drift and constant diffusion, we prove the required
regularized inequality explicitly. All history terms remain evaluated
on the regularized functions.

Under Assumption~\ref{ass:standing}, let $u\in USC(Q)$ and
$v\in LSC(Q)$ be, respectively, a viscosity
subsolution and a viscosity supersolution of \eqref{eq:pde}, satisfying
\[
u(g,t)\le C(1+|g|),\quad v(g,t)\ge-C(1+|g|),
\quad u(g,T)\le0\le v(g,T).
\]
For $\eta>0$, set
\[
M_\eta:=\max_Q\left\{
u(g,t)-v(g,t)-\frac\eta2g^2
\right\}.
\]

\begin{lemma}\label{lem:right-sums}
Suppose $M_\eta>0$ and $(z,t_0)\in Q_0$ is a maximizer satisfying
$B[u](z,t_0)>0$. Then
\begin{equation}\label{eq:right-sums-estimate}
rM_\eta+\eta\left(\rho+\frac r2\right)z^2
\le\frac{\sigma^2}{2}\eta.
\end{equation}
\end{lemma}

\begin{proof}
Fix $a>0$.
Extend each candidate constantly beyond maturity
\[
\check w(g,s):=
\begin{cases}
w(g,s),&0\le s\le T,\\
w(g,T),&s>T,
\end{cases}
\quad w=u,v.
\]
Following \cite[Proposition~3.3]{Namba18} for convolution centers $(g,t)\in\R^2$, we define
\begin{align*}
u^{\varepsilon,\kappa}(g,t)
&:=\sup_{\substack{y\in\R\\s\ge0}}
\left\{\check u(y,s)-\frac{(g-y)^2}{2\varepsilon}
-\frac{(t-s)^2}{2\kappa}\right\},\\
v_{\varepsilon,\kappa}(g,t)
&:=\inf_{\substack{y\in\R\\s\ge0}}
\left\{\check v(y,s)+\frac{(g-y)^2}{2\varepsilon}
+\frac{(t-s)^2}{2\kappa}\right\}.
\end{align*}
For $0<\beta<1$, let
\[
\mathcal D^\beta W(g,t):=
\frac{\beta}{\Gamma(1-\beta)}
\int_0^\infty
\frac{W(g,t)-W(g,t+\tau)}{\tau^{1+\beta}}\,d\tau
\]
whenever the integral is absolutely convergent. This operator is used
on the regularizations; they are not assumed constant after $T$.

We set 
\[
\Theta(t):=
\begin{cases}
(t-t_0)^4,&t\le T,\\
(T-t_0)^4,&t>T.
\end{cases}
\]
The function $\Theta$ is locally Lipschitz on $\R$ and smooth on
$(-\infty,T)$, although its derivative jumps at $T$. For each $A<T$,
$\Theta$ is bounded and Lipschitz on $[A,\infty)$, so
\[
 |\Theta(t+\tau)-\Theta(t)|\le C_A\min\{\tau,1\},
 \qquad t\ge A,\quad \tau>0.
\]
Since $\min\{\tau,1\}\tau^{-1-\beta}\in L^1(0,\infty)$ for $0<\beta<1$,
the integral defining $\mathcal D^\beta\Theta(t)$ converges absolutely
and depends continuously on $t$ by dominated convergence.

Choose positive sequences $\varepsilon_n,\kappa_n,\lambda_n\to0$ and let
$U_n:=u^{\varepsilon_n,\kappa_n}$ and
$V_n:=v_{\varepsilon_n,\kappa_n}$. The function to be maximized is
\begin{align*}
\Psi_n(g,h,t):={}&U_n(g,t)-V_n(h,t)
-\frac{(g-h)^2}{2\lambda_n}-\frac\eta4(g^2+h^2)
-a\bigl((g-z)^4+(h-z)^4+\Theta(t)\bigr).
\end{align*}
We shall show that $\Psi_n$ admits a global maximizer
$(g_n,h_n,t_n)\in\R^3$.

We first record the effect of the extension. For $0\le t<T$, splitting
the integral at $T-t$ gives
\begin{align*}
\mathcal D^\beta\check w(g,t)
&=\frac{w(g,t)-w(g,T)}{\Gamma(1-\beta)(T-t)^\beta}+\frac{\beta}{\Gamma(1-\beta)}
\int_0^{T-t}\frac{w(g,t)-w(g,t+\tau)}{\tau^{1+\beta}}\,d\tau,
\end{align*}
whenever it is finite. Thus this is precisely the right fractional
difference operator for the original candidate.

The growth bounds ensure that the convolutions are finite and their
extrema are attained. For fixed positive $\varepsilon,\kappa$, they
are locally Lipschitz in $(g,t)$ on $\R^2$. They are also bounded
on $[-R,R]\times[A,\infty)$ for every finite $R$ and $A$.
Therefore, for $0<\beta<1$, their $\mathcal D^\beta$ histories
are absolutely convergent and continuous. Indeed, local Lipschitz
continuity bounds the integral over $(0,d)$ by a constant times
$d^{1-\beta}$, while boundedness on the whole future strip gives
dominated convergence on $(d,\infty)$.

Let $t_-:=\max\{-t,0\}$. Since the convolution is taken over
$s\ge0$, we have $|t-s|\ge t_-$. Moreover, $|y|\le |g|+|y-g|$,
and completing the square gives
\[
C|y-g|-\frac{|y-g|^2}{2\varepsilon}
=\frac{C^2\varepsilon}{2}
 -\frac{(|y-g|-C\varepsilon)^2}{2\varepsilon}
\le \frac{C^2\varepsilon}{2}.
\]
Thus, using $\check u(y,s)\le C(1+|y|)$, we obtain
\[
\check u(y,s)-\frac{|g-y|^2}{2\varepsilon}
             -\frac{|t-s|^2}{2\kappa}
\le C(1+|g|)+\frac{C^2\varepsilon}{2}
             -\frac{t_-^2}{2\kappa}.
\]
Taking the supremum over $(y,s)$ gives the upper bound for
$u^{\varepsilon,\kappa}$. Applying the same calculation to
$\check v(y,s)\ge -C(1+|y|)$ gives the lower bound for
$v_{\varepsilon,\kappa}$. Hence
\begin{align*}
u^{\varepsilon,\kappa}(g,t)
&\le C(1+|g|)+\frac{C^2\varepsilon}{2}
             -\frac{t_-^2}{2\kappa},\\
v_{\varepsilon,\kappa}(g,t)
&\ge-C(1+|g|)-\frac{C^2\varepsilon}{2}
             +\frac{t_-^2}{2\kappa}.
\end{align*}

We next localize the positive values of $\Psi_n$ before asserting
that a maximizer exists. Taking $(y,s)=(z,t_0)$ in both convolutions
gives $U_n(z,t_0)\ge u(z,t_0)$ and $V_n(z,t_0)\le v(z,t_0)$. Thus
\[
\Psi_n(z,z,t_0)
\ge u(z,t_0)-v(z,t_0)-\frac{\eta}{2}z^2
=M_\eta>0.
\]
Fix $\varepsilon_0>0$ such that $\varepsilon_n\le\varepsilon_0$
for every $n$. Dropping the nonnegative doubling and quartic
penalties from the definition of $\Psi_n$, we obtain
\begin{align*}
\Psi_n(g,h,t)
&\le 2C+C(|g|+|h|)+C^2\varepsilon_0
     -\frac{\eta}{4}(g^2+h^2)-\frac{t_-^2}{\kappa_n}\\
&\le K_\eta-\frac{\eta}{8}(g^2+h^2)-\frac{t_-^2}{\kappa_n},
\end{align*}
where we used $C|x|\le\eta x^2/8+2C^2/\eta$ and set
$K_\eta:=2C+C^2\varepsilon_0+4C^2/\eta$. In particular,
\[
\Psi_n(g,h,t)>0
\quad\Longrightarrow\quad
|g|,|h|\le R:=\sqrt{8K_\eta/\eta},
\quad t\ge-\sqrt{K_\eta\kappa_n}.
\]

It remains to bound $t$ from above. Set
\[
L_R:=C(1+R)+\frac{C^2\varepsilon_0}{2}.
\]
Fix $|g|\le R$, $t\ge T$, and $0<\varepsilon\le\varepsilon_0$.
If $s\ge T$, the constant extension and the terminal signs give
$\check u(y,s)=u(y,T)\le0$ and $\check v(y,s)=v(y,T)\ge0$.
Thus these source points contribute at most zero to the supremum
and at least zero to the infimum. If $0\le s<T$, then
$|t-s|\ge t-T$, and the preceding spatial estimate gives
\begin{align*}
\check u(y,s)-\frac{|g-y|^2}{2\varepsilon}
             -\frac{|t-s|^2}{2\kappa}
&\le L_R-\frac{(t-T)^2}{2\kappa},\\
\check v(y,s)+\frac{|g-y|^2}{2\varepsilon}
             +\frac{|t-s|^2}{2\kappa}
&\ge-L_R+\frac{(t-T)^2}{2\kappa}.
\end{align*}
Combining the two ranges of $s$ yields
\[
u^{\varepsilon,\kappa}(g,t)
\le\max\left\{0,L_R-\frac{(t-T)^2}{2\kappa}\right\},
\quad
v_{\varepsilon,\kappa}(g,t)
\ge-\max\left\{0,L_R-\frac{(t-T)^2}{2\kappa}\right\}.
\]
Therefore, if $t\ge T+\sqrt{2L_R\kappa_n}$ and $|g|,|h|\le R$,
then $U_n(g,t)\le0\le V_n(h,t)$, so $\Psi_n(g,h,t)\le0$.

The set $\{\Psi_n\ge M_\eta\}$ is therefore nonempty and closed,
by continuity of $\Psi_n$, and is contained in the compact set
\[
[-R,R]^2\times
\bigl[-\sqrt{K_\eta\kappa_n},\,T+\sqrt{2L_R\kappa_n}\bigr].
\]
Hence $\Psi_n$ attains its global maximum. Since $\kappa_n\to0$,
all its maximizers lie in a common compact set.

Let $(g_n,h_n,t_n)$ be any such maximizer. As all the penalties
are nonnegative,
\[
U_n(g_n,t_n)-V_n(h_n,t_n)
\ge\Psi_n(g_n,h_n,t_n)\ge M_\eta>0.
\]
Combining this with the one-sided convolution bounds gives
\[
-L_R\le V_n(h_n,t_n)<U_n(g_n,t_n)\le L_R.
\]
Now let $(y_n,s_n)$ realize the supremum defining $U_n(g_n,t_n)$.
By the defining equality,
\begin{align}
\frac{|g_n-y_n|^2}{2\varepsilon_n}
+\frac{|t_n-s_n|^2}{2\kappa_n}
&=\check u(y_n,s_n)-U_n(g_n,t_n)\le C(1+|y_n|)+L_R\nonumber\\
&\le C(1+R)+L_R+C|g_n-y_n|.\label{ineq: upper}
\end{align}
Using Young's inequality in the form
\[
C|g_n-y_n|
\le\frac{|g_n-y_n|^2}{4\varepsilon_n}+C^2\varepsilon_n,
\]
and  applying to \eqref{ineq: upper}, we conclude  that
\[
\frac{|g_n-y_n|^2}{4\varepsilon_n}
+\frac{|t_n-s_n|^2}{2\kappa_n}
\le C(1+R)+L_R+C^2\varepsilon_0=:D_R.
\]
For an infimum realizer $(\widetilde y_n,\widetilde s_n)$,
the corresponding penalty equals
\[
V_n(h_n,t_n)-\check v(\widetilde y_n,\widetilde s_n)
\le L_R+C(1+|\widetilde y_n|).
\]
The same estimate therefore holds with
$(g_n,y_n,s_n)$ replaced by
$(h_n,\widetilde y_n,\widetilde s_n)$. Taking square roots yields
\begin{align*}
|g_n-y_n|+|h_n-\widetilde y_n|
&\le 4\sqrt{D_R}\,\sqrt{\varepsilon_n},\\
|t_n-s_n|+|t_n-\widetilde s_n|
&\le 2\sqrt{2D_R}\,\sqrt{\kappa_n}.
\end{align*}
The constant $D_R$ is independent of $n$ and of the choice
of realizing points.

The doubled spatial penalty implies $g_n-h_n\to0$. Along any
convergent subsequence, write its limit as $(x,x,t)$, with
$t\in[0,T]$. Semicontinuity of the original candidates and the
realizing-point estimates give
\begin{align*}
M_\eta
&\le\liminf_n\Psi_n(g_n,h_n,t_n)\le\limsup_n\Psi_n(g_n,h_n,t_n)\\
&\le u(x,t)-v(x,t)-\frac\eta2x^2
-2a(x-z)^4-a\Theta(t)\\
&\le M_\eta-2a(x-z)^4-a\Theta(t).
\end{align*}
Thus $(x,t)=(z,t_0)$ and the doubled penalty tends to zero. Hence
\[
(g_n,h_n,t_n)\longrightarrow(z,z,t_0),\quad
U_n(g_n,t_n)\longrightarrow u(z,t_0),\quad
V_n(h_n,t_n)\longrightarrow v(z,t_0).
\]
In particular, all realizing points converge to $(z,t_0)$ and
belong to $Q_0$ for large $n$. Their original values converge as
well, for example,
\[
U_n(g_n,t_n)\le u(y_n,s_n),\quad
\limsup_{n\to\infty}u(y_n,s_n)\le u(z,t_0).
\]

Let $b:=B[u](z,t_0)>0$. In a sufficiently small neighborhood of
$(z,t_0)$, upper semicontinuity at the reflected point and
continuity of $\psi$ give
\[
u(-y,s)\le u(-z,t_0)+\frac b4,
\quad
|\psi(y,s)-\psi(z,t_0)|\le\frac b4.
\]
At all large-$n$ supremum realizers, $u(y,s)\ge U_n(g,t)$.
Since $U_n(g_n,t_n)\to u(z,t_0)$, continuity of $U_n$ gives a
neighborhood of $(g_n,t_n)$ in which
$U_n(g,t)>u(z,t_0)-b/4$. Shrinking this neighborhood if necessary,
the realizing-point bounds keep every source point in the preceding
fixed neighborhood. Therefore
\[
\begin{aligned}
B[u](y,s)
&=u(y,s)-u(-y,s)-\psi(y,s)\\
&>\left(u(z,t_0)-\frac b4\right)
 -\left(u(-z,t_0)+\frac b4\right)
 -\left(\psi(z,t_0)+\frac b4\right)\\
&=b-\frac{3b}{4}=\frac b4>0.
\end{aligned}
\]
at every supremum realizer for these centers.

Fix $n$ sufficiently large. Let $P_n^+=(g_n^+,t_n^+)$ be any center
in the neighborhood of $(g_n,t_n)$ obtained above, and let $\varphi_n^+$
touch $U_n=u^{\varepsilon_n,\kappa_n}$ from above at $P_n^+$, with
$\varphi_n^+(P_n^+)=U_n(P_n^+)$. If $(y_n^+,s_n^+)$ realizes the
supremum at $P_n^+$, set
\[
 d_n^+:=\frac{(g_n^+ - y_n^+)^2}{2\varepsilon_n}
       +\frac{(t_n^+ - s_n^+)^2}{2\kappa_n}.
\]
Then $U_n(P_n^+)=u(y_n^+,s_n^+)-d_n^+$, and
\[
 (x,\tau)\longmapsto
 \varphi_n^+(g_n^+ + x-y_n^+,\,t_n^+ + \tau-s_n^+)+d_n^+
\]
is a normalized local upper test for $u$ at $(y_n^+,s_n^+)$.
Keeping this realizer fixed and varying the spatial center gives
\[
 (\varphi_n^+)_g(P_n^+)=\frac{y_n^+ - g_n^+}{\varepsilon_n},
 \qquad y_n^+=g_n^+ + \varepsilon_n(\varphi_n^+)_g(P_n^+).
\]

For $0<\beta<1$, we show that the actual history at $(y_n^+,s_n^+)$
is absolutely convergent. Write
\[
 d\mu_\beta(\tau):=\frac{\beta}{\Gamma(1-\beta)}\frac{d\tau}{\tau^{1+\beta}},
 \qquad
 I_n(\tau):=u(y_n^+,s_n^+)-\check u(y_n^+,s_n^+ + \tau).
\]
Choose $0<\delta_n<T-s_n^+$ so that the translated upper test
dominates $u$ at $(y_n^+,s_n^+ + \tau)$ for $0\le \tau\le\delta_n$.
At $\tau=0$ there is equality. Subtracting the upper bound
$u(y_n^+,s_n^+ + \tau)\le\varphi_n^+(g_n^+,t_n^+ + \tau)+d_n^+$
from $u(y_n^+,s_n^+)=\varphi_n^+(P_n^+)+d_n^+$ gives
\begin{align*}
 I_n(\tau)&\ge J_n(\tau):=\varphi_n^+(P_n^+)
                     -\varphi_n^+(g_n^+,t_n^+ + \tau),\\
 J_n(\tau)&=-\int_0^\tau(\varphi_n^+)_t(g_n^+,t_n^+ + q)\,dq,
 \qquad |J_n(\tau)|\le C_n \tau,
\end{align*}
where $C_n:=\max_{0\le q\le\delta_n}
|(\varphi_n^+)_t(g_n^+,t_n^+ + q)|<\infty$.
Writing $I_n^\pm(\tau):=\max\{\pm I_n(\tau),0\}$, we obtain, for
$0<\delta\le\delta_n$,
\begin{equation}\label{eq:right-sums-near-test}
 \int_0^\delta I_n^-\,d\mu_\beta
 \le\int_0^\delta |J_n|\,d\mu_\beta
 \le\frac{\beta C_n}{\Gamma(1-\beta)}
      \int_0^\delta \tau^{-\beta}\,d\tau
 =\frac{\beta C_n}{(1-\beta)\Gamma(1-\beta)}\,\delta^{1-\beta}.
\end{equation}
For $\tau\ge\delta_n$, the growth bound
$\check u(y_n^+,s_n^+ + \tau)\le C(1+|y_n^+|)$ bounds $I_n^-(\tau)$
from above. Since $\int_{\delta_n}^{\infty}d\mu_\beta<\infty$, this tail
bound and \eqref{eq:right-sums-near-test} with $\delta=\delta_n$ give
$N_n:=\int_0^\infty I_n^-\,d\mu_\beta<\infty$.

To control the positive part, put
\[
 \Lambda_n:=\rho y_n^+(\varphi_n^+)_g(P_n^+)
 -\frac{\sigma^2}{2}(\varphi_n^+)_{gg}(P_n^+)+r u(y_n^+,s_n^+).
\]
Since $B[u](y_n^+,s_n^+)>0$, the differential branch of the split
subsolution inequality holds for every $0<\delta<\delta_n,$ i.e.
\begin{equation}\label{eq:right-sums-split-upper}
 \int_0^\delta J_n\,d\mu_\beta
 +\int_\delta^\infty I_n\,d\mu_\beta+\Lambda_n\le0.
\end{equation}
Here the constant extension $\check u$ includes the terminal term
\[
 \int_{T-s_n^+}^\infty I_n\,d\mu_\beta
 =\frac{u(y_n^+,s_n^+)-u(y_n^+,T)}
        {\Gamma(1-\beta)(T-s_n^+)^\beta}.
\]
Substituting $I_n=I_n^+-I_n^-$ in \eqref{eq:right-sums-split-upper}
and applying \eqref{eq:right-sums-near-test} to the test integral yields
\begin{equation}\label{eq:right-sums-positive-part}
\begin{aligned}
 \int_\delta^\infty I_n^+\,d\mu_\beta
 &\le -\Lambda_n-\int_0^\delta J_n\,d\mu_\beta
       +\int_\delta^\infty I_n^-\,d\mu_\beta\\
 &\le |\Lambda_n|+N_n
       +\frac{\beta C_n}{(1-\beta)\Gamma(1-\beta)}
         \,\delta_n^{1-\beta}.
\end{aligned}
\end{equation}
The bound in \eqref{eq:right-sums-positive-part} is independent of
$\delta$. Letting $\delta\downarrow0$ there and applying monotone
convergence gives $\int_0^\infty I_n^+\,d\mu_\beta<\infty$.
Together with $N_n<\infty$, this gives
$\int_0^\infty |I_n|\,d\mu_\beta<\infty$.
We may now pass to the limit in \eqref{eq:right-sums-split-upper}.
The test integral tends to zero by \eqref{eq:right-sums-near-test},
and the truncated integral tends to
$\int_0^\infty I_n\,d\mu_\beta=\mathcal D^\beta\check u(y_n^+,s_n^+)$.
Hence
\begin{equation}\label{eq:right-sums-source-upper}
 \mathcal D^\beta\check u(y_n^+,s_n^+)
 +\rho y_n^+(\varphi_n^+)_g(P_n^+)
 -\frac{\sigma^2}{2}(\varphi_n^+)_{gg}(P_n^+)
 +r u(y_n^+,s_n^+)\le0.
\end{equation}
For every $\tau\ge0$, using $(y_n^+,s_n^+ + \tau)$ as a competitor at
$(g_n^+,t_n^+ + \tau)$ leaves the penalty $d_n^+$ unchanged. Therefore
\begin{align*}
 U_n(g_n^+,t_n^+ + \tau)
 &\ge \check u(y_n^+,s_n^+ + \tau)-d_n^+,\\
 U_n(P_n^+)-U_n(g_n^+,t_n^+ + \tau)
 &\le \check u(y_n^+,s_n^+)-\check u(y_n^+,s_n^+ + \tau).
\end{align*}
Integrating gives
$\mathcal D^\beta U_n(P_n^+)\le \mathcal D^\beta\check u(y_n^+,s_n^+)$.
Combining this comparison with \eqref{eq:right-sums-source-upper}
and using the contact value and gradient identity above, we obtain
\begin{equation}\label{ineq:super}
\begin{gathered}
     \mathcal D^\beta U_n(P_n^+)+\rho g_n^+(\varphi_n^+)_g(P_n^+)
 -\frac{\sigma^2}{2}(\varphi_n^+)_{gg}(P_n^+)+rU_n(P_n^+)\\
 \qquad\le -\rho\varepsilon_n\lvert(\varphi_n^+)_g(P_n^+)\rvert^2
             -rd_n^+\le0.
\end{gathered}
\end{equation}
For $\beta=1$, the translated test has the same time derivative at the
source, so the last inequality holds with $-(\varphi_n^+)_t(P_n^+)$
in place of $\mathcal D^\beta U_n(P_n^+)$.

Similarly, let $P_n^-=(h_n^-,t_n^-)$ be any center in the neighborhood
of $(h_n,t_n)$ obtained above, and let $\varphi_n^-$ be a normalized
local lower test for $V_n=v_{\varepsilon_n,\kappa_n}$ at $P_n^-$.
If $(y_n^-,s_n^-)$ realizes the infimum at $P_n^-$, set
\[
 d_n^-:=\frac{(h_n^- - y_n^-)^2}{2\varepsilon_n}
       +\frac{(t_n^- - s_n^-)^2}{2\kappa_n}.
\]
The translated lower test is
\[
 (x,\tau)\longmapsto
 \varphi_n^-(h_n^- + x-y_n^-,\,t_n^- + \tau-s_n^-)-d_n^-,
\]
and the corresponding identities are
\[
 y_n^-=h_n^- - \varepsilon_n(\varphi_n^-)_g(P_n^-),
 \qquad V_n(P_n^-)=v(y_n^-,s_n^-)+d_n^-.
\]
For $0<\beta<1$, the lower test makes the positive part of the history
increment integrable near zero, and the split supersolution inequality
bounds the negative part. Using $(y_n^-,s_n^- + \tau)$ as a competitor at
$(h_n^-,t_n^- + \tau)$ now gives
$\mathcal D^\beta V_n(P_n^-)\ge\mathcal D^\beta\check v(y_n^-,s_n^-)$.
Therefore, we get
\begin{equation}\label{ineq:sub}   
\begin{gathered}
 \mathcal D^\beta V_n(P_n^-)+\rho h_n^-(\varphi_n^-)_g(P_n^-)
 -\frac{\sigma^2}{2}(\varphi_n^-)_{gg}(P_n^-)+rV_n(P_n^-)\\
 \qquad\ge \rho\varepsilon_n\lvert(\varphi_n^-)_g(P_n^-)\rvert^2
             +rd_n^-\ge0.
\end{gathered}
\end{equation}   
For $\beta=1$, the same inequality holds with $-(\varphi_n^-)_t(P_n^-)$
in place of $\mathcal D^\beta V_n(P_n^-)$.

The convolution centers range over $\R^2$. Therefore, the translated
test is defined on a full neighborhood of its center even if $t<0$;
at an original source with $s=0$, it is a valid relative test on $Q$.
This proves the regularized inequalities also when $t_0=0$.

Since $t_n\to t_0<T$, a fixed open interval $I$ with
$t_0\in I$ and $\sup I<T$ contains $t_n$ for all sufficiently large $n$.
On $I$, $\Theta(t)=(t-t_0)^4$, so the penalty in $\Psi_n$ is smooth.
We apply the parabolic theorem of sums \cite[Theorem~8.3]{CIL92}
to $U_n$ and $-V_n$ at $(g_n,h_n,t_n)$ on this interval.

Local time Lipschitz continuity supplies the
required bounds on temporal jet components. For $0<\beta<1$,
continuity of the regularized histories permits the differential
inequalities just proved to pass to closed parabolic jets.
Thus, for all sufficiently large $n$, there are closed parabolic jets
\[
(\ell_n,p_n,X_n)\in\overline{\mathcal P}^{2,+}U_n(g_n,t_n),
\quad
(m_n,q_n,Y_n)\in\overline{\mathcal P}^{2,-}V_n(h_n,t_n),
\]
where the first component is the temporal slope, such that
\[
\ell_n-m_n=a\Theta'(t_n),
\]
and
\[
\operatorname{diag}(X_n,-Y_n)\le A_n+\lambda_n A_n^2,
\]
where
\[
A_n=\frac1{\lambda_n}
\begin{pmatrix}1&-1\\-1&1\end{pmatrix}
+\begin{pmatrix}
\frac\eta2+12a(g_n-z)^2&0\\
0&\frac\eta2+12a(h_n-z)^2
\end{pmatrix}.
\]
The gradients are
\[
p_n=\frac{g_n-h_n}{\lambda_n}+\frac\eta2g_n+4a(g_n-z)^3,
\quad
q_n=\frac{g_n-h_n}{\lambda_n}-\frac\eta2h_n-4a(h_n-z)^3.
\]
Define the history evaluation by
\[
\mathcal H_n:=
\begin{cases}
\mathcal D^\beta U_n(g_n,t_n)-\mathcal D^\beta V_n(h_n,t_n),
&0<\beta<1,\\
-\ell_n+m_n,&\beta=1.
\end{cases}
\]

Fix $n$ sufficiently large. The regularized inequalities \eqref{ineq:super} and \eqref{ineq:sub} hold for ordinary parabolic superjets
of $U_n$ near $(g_n,t_n)$ and subjets of $V_n$ near $(h_n,t_n)$,
by the test-function characterization of these jets. By definition of
the closed jets selected above, there are such ordinary jets whose base
points and components satisfy, as $k\to\infty$,
\[
\begin{aligned}
 (g_{n,k}^+,t_{n,k}^+;\ell_{n,k},p_{n,k},X_{n,k})
 &\longrightarrow (g_n,t_n;\ell_n,p_n,X_n),\\
 (h_{n,k}^-,t_{n,k}^-;m_{n,k},q_{n,k},Y_{n,k})
 &\longrightarrow (h_n,t_n;m_n,q_n,Y_n).
\end{aligned}
\]
The drift and discount corrections are nonpositive for $U_n$ and
nonnegative for $V_n$, so they may be dropped before taking the limit.
For $0<\beta<1$, the functions $U_n,V_n$ and their histories
$\mathcal D^\beta U_n,\mathcal D^\beta V_n$ are continuous for this
fixed $n$. Passing to the limit as $k\to\infty$ therefore gives
\[
\begin{aligned}
 \mathcal D^\beta U_n(g_n,t_n)
 +\rho g_np_n-\frac{\sigma^2}{2}X_n+rU_n(g_n,t_n)&\le0,\\
 \mathcal D^\beta V_n(h_n,t_n)
 +\rho h_nq_n-\frac{\sigma^2}{2}Y_n+rV_n(h_n,t_n)&\ge0.
\end{aligned}
\]
For $\beta=1$, the same inequalities hold with time terms $-\ell_n$
and $-m_n$, respectively, since
$-\ell_{n,k}\to-\ell_n$ and $-m_{n,k}\to-m_n$.
In either case the first left-hand side is at most the second.
Subtracting the second from the first and using the definition of
$\mathcal H_n$ yields
\begin{equation}\label{eq:right-sums-ineq}
 \mathcal H_n+\rho(g_np_n-h_nq_n)
 -\frac{\sigma^2}{2}(X_n-Y_n)
 +r\bigl(U_n(g_n,t_n)-V_n(h_n,t_n)\bigr)\le0.
\end{equation}

Set
\[
K_\beta:=
\begin{cases}
\displaystyle\frac{4(T-t_0)^{4-\beta}}
{(4-\beta)\Gamma(1-\beta)},&0<\beta<1,\\[1ex]
0,&\beta=1.
\end{cases}
\]
For the fractional history estimate, set
$W_n(s):=U_n(g_n,s)-V_n(h_n,s)$.
Global maximality in the center time gives, for every $\tau\ge0$,
\[
W_n(t_n)-W_n(t_n+\tau)
\ge a\bigl(\Theta(t_n)-\Theta(t_n+\tau)\bigr).
\]
All these histories are absolutely convergent. Therefore
\[
\mathcal H_n\ge a\mathcal D^\beta\Theta(t_n)
\longrightarrow a\mathcal D^\beta\Theta(t_0)=-aK_\beta.
\]
Indeed, with $L:=T-t_0$,
\[
\mathcal D^\beta\Theta(t_0)
=-\frac{\beta}{\Gamma(1-\beta)}
\left(\frac{L^{4-\beta}}{4-\beta}
+\frac{L^{4-\beta}}{\beta}\right)=-K_\beta.
\]
For $\beta=1$,
$\mathcal H_n=-a\Theta'(t_n)\to0$.
Therefore, for every $\beta\in(0,1]$,
\begin{equation}\label{eq:right-sums-history}
\liminf_{n\to\infty}\mathcal H_n\ge-aK_\beta.
\end{equation}
Take $e=(1,1)^\top$, then the singular doubled part of $A_n$ annihilates
$e$, so $A_ne$ is bounded and
\[
X_n-Y_n\le e^\top A_ne+\lambda_n|A_ne|^2
=\eta+o(1).
\]
Also
\[
g_np_n-h_nq_n
=\frac{(g_n-h_n)^2}{\lambda_n}
+\frac\eta2(g_n^2+h_n^2)
+4a\bigl(g_n(g_n-z)^3+h_n(h_n-z)^3\bigr).
\]
Taking limits we obtain
\begin{align}
\limsup_{n\to\infty}(X_n-Y_n)&\le\eta,
\label{eq:right-sums-matrix}\\
\liminf_{n\to\infty}(g_np_n-h_nq_n)&\ge\eta z^2,
\label{eq:right-sums-drift}\\
U_n(g_n,t_n)-V_n(h_n,t_n)&\longrightarrow
M_\eta+\frac\eta2z^2.
\label{eq:right-sums-value}
\end{align}
Using these bounds and \eqref{eq:right-sums-history} in
\eqref{eq:right-sums-ineq}, with $\rho\ge0$, yields
\[
0\ge
-aK_\beta+\rho\eta z^2-\frac{\sigma^2}{2}\eta
+r\left(M_\eta+\frac\eta2z^2\right).
\]
Since $a>0$ is arbitrary and $K_\beta$ is finite and independent of $a$,
letting $a\to0^+$ proves \eqref{eq:right-sums-estimate}.
\end{proof}

\begin{lemma}
\label{lem:selection-sc}
Let $u\in USC(Q)$ and $v\in LSC(Q)$ satisfy $B[v]\ge0$ on $Q_0$. If
$u-v-\eta g^2/2$ has a global maximum at $(x,t)$ with $t<T$, then either
$B[u](x,t)>0$, or $(-x,t)$ is also a global maximum and
$B[u](-x,t)=2c>0$.
\end{lemma}

\begin{proof}
If $B[u](x,t)>0$, there is nothing to prove. Otherwise,
\[
(u-v)(x,t)-(u-v)(-x,t)=B[u](x,t)-B[v](x,t)\le0.
\]
The penalty is even, so maximality forces equality and makes $(-x,t)$
another maximizer. Since $B[u](x,t)\le0\le B[v](x,t)$, equality implies
$B[u](x,t)=B[v](x,t)=0$. The identity
$B[u](x,t)+B[u](-x,t)=2c$ then gives $B[u](-x,t)=2c$.
\end{proof}

\begin{theorem}[Comparison for semicontinuous solutions]
\label{thm:comp}
Under Assumption~\ref{ass:standing}, let $u\in USC(Q)$ be a viscosity
subsolution and let $v\in LSC(Q)$ be a viscosity supersolution of
\eqref{eq:pde}. Assume the uniform linear-growth bounds
\[
u(g,t)\le C(1+|g|),\quad v(g,t)\ge-C(1+|g|),
\]
and the terminal inequalities $u(g,T)\le0\le v(g,T)$. Then
\[
u\le v\quad\text{on }Q.
\]
\end{theorem}

\begin{proof}
Assume that $u-v$ is positive somewhere. Choose
$(g_*,t_*)$ and $m>0$ such that
$(u-v)(g_*,t_*)=m$. For $\eta>0$ sufficiently small,
\[
M_\eta:=\max_Q\left\{u(g,t)-v(g,t)-\frac\eta2g^2\right\}>\frac m2.
\]
The maximum is finite and attained because of the linear-growth bounds.
It cannot occur at $t=T$, so let $(z,t_0)$ be a maximizer with $t_0<T$.
By Lemma~\ref{lem:selection-sc}, after replacing $z$ by $-z$ if necessary,
we may assume
\[
B[u](z,t_0)>0.
\]

By \eqref{eq:right-sums-estimate}, we have $rM_\eta\le \frac{\sigma^2}{2}\eta$.
For $0<\eta<rm/\sigma^2$, this contradicts $M_\eta>m/2$.
Hence 
$$
u\le v\quad\text{on}\quad Q.
$$
\end{proof}

\begin{corollary}[Uniqueness]
\label{cor:unique}
There is at most one viscosity solution of \eqref{eq:pde} in the uniform
linear-growth class. In particular, any two continuous viscosity
solutions coincide.
\end{corollary}

\begin{proof}
Apply Theorem~\ref{thm:comp} first with the subsolution envelope of the
first solution and the supersolution envelope of the second one, and then
with the roles reversed.
\end{proof}

\section{Existence of Viscosity Solutions}\label{sec:exist}

We examine Perron's method for the existence question. The barriers below
provide the lower and upper bounds, and the semicontinuous comparison theorem
above supplies the final envelope step.

\subsection{The subsolution barrier}\label{ssec:sub}

\begin{lemma}\label{lem:sub}
The non-negative floor of the obstacle profile, defined by 
$$
\underline u(g,t):=\max(0,\psi(g,t)),
$$ 
serves as a viscosity subsolution of the fractional bubble problem \eqref{eq:pde} on $\R \times [0,T)$.
\end{lemma}

\begin{proof}
To establish the viscosity subsolution property, we choose an arbitrary test point $(g_0,t_0)$ with $t_0<T$. Let $\varphi \in C^{2,1}(\R \times [0,T])$ be a smooth test function such that $\underline u - \varphi$ attains a local maximum at $(g_0,t_0)$. Without loss of generality, we normalize the function such that $\underline u(g_0,t_0) = \varphi(g_0,t_0)$, implying $\underline u(g,t) \le \varphi(g,t)$ in the contact neighborhood. 

We must demonstrate that the subsolution inequality holds at the contact point
\begin{equation}\label{eq:sub_def}
    \min\Bigl\{F_{\varphi,\delta}(g_0,t_0;\underline u),\;\; \underline u(g_0,t_0) - \underline u(-g_0,t_0) - \psi(g_0,t_0)\Bigr\} \le 0.
\end{equation}
We partition the analysis into two mutually exclusive regimes based on the sign of the obstacle profile at the contact point.

\textbf{\emph{Case 1: ($\psi(g_0,t_0) \ge 0$).}}
By definition, the value function binds exactly to the obstacle, yielding $\underline u(g_0,t_0) = \psi(g_0,t_0)$. We substitute this into the obstacle branch of the minimum operator in \eqref{eq:sub_def} to evaluate the obstacle condition directly
\begin{align*}
    \underline u(g_0,t_0) - \underline u(-g_0,t_0) - \psi(g_0,t_0) 
    &= \psi(g_0,t_0) - \max\bigl(0, \psi(-g_0,t_0)\bigr) - \psi(g_0,t_0) \\
    &= -\max\bigl(0, \psi(-g_0,t_0)\bigr) \le 0.
\end{align*}
Since the obstacle branch evaluates to a non-positive value, the minimum operator in \eqref{eq:sub_def} is automatically $\le 0$. The subsolution property is verified in this regime without requiring evaluation of the differential operator $F_{\varphi,\delta}$.

\textbf{\emph{Case 2: ($\psi(g_0,t_0) < 0$).}}
 Here, the value function satisfies $\underline u(g_0,t_0) = 0$. Due to $\underline u(g,t) \ge 0$ everywhere by construction, the smooth test function satisfies, in the contact neighborhood,
\[
    \varphi(g,t) \ge \underline u(g,t) \ge 0 = \varphi(g_0,t_0).
\]
Therefore, the contact point $(g_0,t_0)$ transforms into a local spatial minimum for the smooth function $\varphi(\cdot,t_0)$. From classical first- and second-order optimality conditions, the spatial derivatives must satisfy
\begin{equation}\label{eq:spatial_opt}
    \varphi_g(g_0,t_0) = 0, \quad \text{and} \quad \varphi_{gg}(g_0,t_0) \ge 0.
\end{equation}

Next, we calculate the temporal trajectory along the fixed coordinate $g_0$. The time profile of the obstacle is governed by $\psi(g_0,s) = g_0\alpha(s) - c$. Since $\alpha(s)$ is non-increasing with $\alpha(T)=0$, we analyze the path components
\begin{itemize}
    \item For $g_0 \ge 0$, the mapping $s \mapsto \psi(g_0,s)$ is monotonically non-increasing, which forces $\psi(g_0,s) \le \psi(g_0,t_0) < 0$ for all $s \in [t_0, T]$.
    \item For $g_0 < 0$, the mapping $s \mapsto \psi(g_0,s)$ is monotonically non-decreasing, terminating at $\psi(g_0,T) = -c < 0$. This forces $\psi(g_0,s) < 0$ for all $s \in [t_0, T]$.
\end{itemize}
Thus, $\psi(g_0,s) < 0$ uniformly on the future interval $[t_0,T]$, which mathematically fixes the value function to $\underline u(g_0,s) \equiv 0$ on this segment.

For $0<\beta<1$, evaluate the split history operator. Its terminal
and distant-time terms vanish because $\underline u(g_0,s)=0$ throughout
$[t_0,T]$. At a local upper contact, $\varphi(g_0,t_0)=0$ and
$\varphi(g_0,t_0+\tau)\ge0$ on each admissible near-time segment. Hence
\[
D_\delta[\underline u,\varphi](g_0,t_0)
=-\frac{\beta}{\Gamma(1-\beta)}\int_0^\delta
\frac{\varphi(g_0,t_0+\tau)}{\tau^{1+\beta}}\,d\tau\le0.
\]
For $\beta=1$, an interior contact gives $\varphi_t(g_0,t_0)=0$.
At $t_0=0$, the relative upper contact gives
$\varphi_t(g_0,0)\ge0$. Thus $-\varphi_t(g_0,t_0)\le0$
in both cases.

Finally, we assemble the complete backward differential operator $F_{\varphi,\delta}(g_0,t_0;\underline u)$. Combining our derivative bounds and taking into account \eqref{eq:spatial_opt} we get
\begin{align*}
    F_{\varphi,\delta}(g_0,t_0;\underline{u}) &= D_\delta[\underline u,\varphi](g_0,t_0) + \rho g_0\varphi_g(g_0,t_0) - \frac{\sigma^2}{2}\varphi_{gg}(g_0,t_0) + r \underline{u}(g_0,t_0) \\
    &= D_\delta[\underline u,\varphi](g_0,t_0)  - \frac{\sigma^2}{2}\varphi_{gg}(g_0,t_0) \le 0.  
\end{align*}
The differential branch is non-positive, satisfying the minimum condition \eqref{eq:sub_def}. 

At the terminal boundary $t=T$, the expression satisfies $\underline u(g,T) = \max(0, \psi(g,T)) = \max(0, -c) = 0$, matching the boundary criterion. Thus, $\underline u$ is globally a viscosity subsolution.
\end{proof}

\subsection{The supersolution barrier}\label{ssec:super}

\begin{lemma}\label{lem:super}
Fix any $M_0>0$ and define the spatial profile $f(g):=\sqrt{g^2+M_0}+\tfrac{g}{2}+K$ for some constant $K>0$. Then, for $K$ sufficiently large, the function $\tilde u(g,t):=\alpha(t)f(g)$ is a global viscosity supersolution of the fractional bubble problem \eqref{eq:pde} residing in the linear growth class.
\end{lemma}

\begin{proof}
Let $h(g) := \sqrt{g^2+M_0} + \frac{g}{2}$, such that $f(g) = h(g) + K$. Note that $h(g)$ is strictly positive with a global minimum of $\frac{\sqrt{3}}{2}\sqrt{M_0}$, ensuring $f(g) \ge K > 0$. Furthermore, $f'(g) = h'(g) = \frac{g}{\sqrt{g^2+M_0}} + \frac{1}{2}$, yielding the uniform global bound $\|f'\|_\infty = \|h'\|_\infty=\frac{3}{2}$. The second derivative evaluates to $f''(g) = h''(g) = M_0(g^2+M_0)^{-3/2}$, which attains its global maximum at the origin, strictly bounding the convexity, i.e.
$$
0\le h''(g)\le \frac{1}{\sqrt{M_0}}.
$$
To be a global supersolution, $\tilde u$ must dominate the obstacle everywhere. Evaluating the spatial difference leverages the asymmetric linear term to perfectly cancel the state variable:
\[
    h(g) - h(-g) = \Bigl(\sqrt{g^2+M_0} + \frac{g}{2}\Bigr) - \Bigl(\sqrt{g^2+M_0} - \frac{g}{2}\Bigr) = g.
\]
Since the constant $K$ drops out of the difference, we have exactly $\tilde u(g,t) - \tilde u(-g,t) = \alpha(t)g$. The obstacle is defined as $\psi(g,t) = \alpha(t)g - c$. Therefore, the gap evaluates to
\[
    \tilde u(g,t) - \tilde u(-g,t) - \psi(g,t) = \alpha(t)g - (\alpha(t)g - c) = c.
\]
Due to the transaction cost $c > 0$, the supersolution strictly dominates the obstacle across the entire domain. The obstacle branch of the minimum operator is thus strictly positive, forcing us to satisfy the differential inequality $F[\tilde u] \ge 0$.

Since the spatial operator $\mathcal{A}$ acts only on $g$ and the time operator $\rCder{\beta}_t$ acts only on $t$, we can explicitly compute the linear pricing operator $F[\tilde u]$ and obtain
\[
    F[\tilde u] = \bigl(\rCder{\beta}_{t}\alpha(t)\bigr)f(g) + \alpha(t)\Bigl[\rho g f'(g) - \frac{\sigma^2}{2}f''(g) + r f(g)\Bigr].
\]
To manage the drift term $\rho g f'(g) = \rho g h'(g)$, we use algebraic substitution. By writing $g^2 = (g^2+M_0) - M_0$, we isolate $h(g)$ as follows
\begin{multline*}
       g h'(g) = \frac{g^2}{\sqrt{g^2+M_0}} + \frac{g}{2} = \frac{g^2+M_0-M_0}{\sqrt{g^2+M_0}} + \frac{g}{2} =\\= \Bigl(\sqrt{g^2+M_0} + \frac{g}{2}\Bigr) - \frac{M_0}{\sqrt{g^2+M_0}} = h(g) - \frac{M_0}{\sqrt{g^2+M_0}}.
\end{multline*}
Substituting this drift identity and $f(g) = h(g) + K$ into the bracketed operator terms yields
\[
    F[\tilde u] = \bigl(\rCder{\beta}_{t}\alpha(t)\bigr)f(g) + \alpha(t)\left[ (r+\rho)h(g) + rK - \frac{\rho M_0}{\sqrt{g^2+M_0}} - \frac{\sigma^2}{2}h''(g) \right].
\]
We now establish the positivity of this evaluation. For the first term, we know  $\rCder{\beta}_{t}\alpha(t) \ge 0$, and we established $f(g) \ge \frac{\sqrt{3M_0}}{2} + K > 0$, making the product  non-negative.

For the second term, since the bubble decay profile satisfies $\alpha(t) \ge 0$, it suffices to make the bracketed expression non-negative. Due to $r,\rho \ge 0$ and $h(g) > 0$, the term $(r+\rho)h(g)$ is strictly positive. We bound the remaining negative elements using the structural limits derived above
\[
    \frac{\rho M_0}{\sqrt{g^2+M_0}} \le \frac{\rho M_0}{\sqrt{M_0}} = \rho\sqrt{M_0}, \quad \text{and} \quad \frac{\sigma^2}{2}h''(g) \le \frac{\sigma^2}{2\sqrt{M_0}}.
\]
To guarantee the bracket is globally non-negative, the constant $K$ must absorb these maximum bounds. We simply demand
\[
    rK \ge \rho\sqrt{M_0} + \frac{\sigma^2}{2\sqrt{M_0}} \implies K \ge \frac{1}{r}\left(\rho\sqrt{M_0} + \frac{\sigma^2}{2\sqrt{M_0}}\right).
\]
With $K$ satisfying this threshold, the bracket will be $\ge 0$, forcing $F[\tilde u] \ge 0$.

By the definition of the speculative time profile, $\alpha(T) = 0$. Therefore, at maturity, the candidate satisfies $\tilde u(g,T) = 0 \cdot f(g) = 0$. This trivially satisfies the required supersolution terminal boundary condition $\tilde u(\cdot, T) \ge 0$, completing the proof.
For $\beta=1$, the sign used above is $-\alpha'\ge0$.
At a normalized local lower contact with this smooth candidate, the near-time
test increments are no smaller than those of $\tilde u$, the spatial
gradients agree, and $\varphi_{gg}\le\tilde u_{gg}$. Thus
$F_{\varphi,\delta}(g_0,t_0;\tilde u)\ge F[\tilde u](g_0,t_0)\ge0$,
which verifies the split viscosity condition.
\end{proof}

\begin{remark}
Memory enters only through the single factor $\rCder{\beta}_t\alpha\ge0$, which
helps the supersolution property to hold.  No additional parameter restriction beyond the standing assumptions is needed, for every $\beta\in(0,1]$.
\end{remark}

\subsection{Existence via Perron's method}\label{ssec:perron}

Throughout this subsection we take global tests
from
\[
\mathcal T:=C(Q)\cap C^{2,1}(Q_0).
\]
For $0<\beta<1$, their right fractional operator is understood through the
Marchaud formula \eqref{eq:marchaud}; for $\beta=1$ it is $-\partial_t$.
Time derivatives at $t=0$ are one-sided. The linear-growth condition is
imposed on subsolutions and supersolutions, not on test functions.
Semicontinuous envelopes are taken relative to $Q$.

\begin{theorem}[Existence via Perron's method]\label{thm:exist}
Under Assumption~\ref{ass:standing} and the assumptions on $\alpha$ in
Section~\ref{ssec:struct}, there exists a unique continuous viscosity
solution $u$ of \eqref{eq:pde} with uniform at most linear growth such that
\[
\underline u\le u\le\tilde u
\quad\text{on }\R\times[0,T],
\quad u(g,T)=0\quad(g\in\R).
\]
\end{theorem}

\begin{proof}
Set
\[
B[w](g,t):=w(g,t)-w(-g,t)-\psi(g,t).
\]
We use the notation $F_\varphi$ from \eqref{eq:testop}; its fractional
term is evaluated on the globally touching smooth test function.

If $\psi(g,t)\le0$, then
\[
\underline u(g,t)=0\le\alpha(t)f(g)=\tilde u(g,t).
\]
If $\psi(g,t)>0$, Lemma~\ref{lem:super} gives
\[
\tilde u(g,t)
=\tilde u(-g,t)+\psi(g,t)+c
\ge\psi(g,t)+c>\underline u(g,t).
\]
Thus $0\le\underline u\le\tilde u$ on $Q$. Since
$\alpha(T)=0$ and $\alpha\in C^1([0,T])$, for each $R>0$,
\begin{equation}\label{eq:perron_terminal_bound}
0\le\underline u(g,t)\le\tilde u(g,t)
\le C_R(T-t),\quad |g|\le R,
\quad
C_R:=\|\alpha'\|_\infty\max_{|x|\le R}f(x).
\end{equation}

Define
\[
\mathcal S:=\left\{
 w\in\operatorname{USC}(Q):
 \begin{array}{l}
 w\text{ is a viscosity subsolution of }\eqref{eq:pde},\\
 \underline u\le w\le\tilde u
 \end{array}
 \right\},
\quad
u(g,t):=\sup_{w\in\mathcal S}w(g,t).
\]
Lemma~\ref{lem:sub} yields $\underline u\in\mathcal S$, so
$\mathcal S\ne\varnothing$. Continuity of the barriers implies
\begin{equation}\label{eq:perron_envelope_bounds}
\underline u\le u_*\le u\le u^*\le\tilde u.
\end{equation}

We verify the upper-envelope stability step explicitly. Let $\chi\in\mathcal T$
touch $u^*$ globally from above at $(g_0,t_0)\in Q$, with $t_0<T$ and
$\chi(g_0,t_0)=u^*(g_0,t_0)$. For fixed $\varepsilon>0$, set
\[
P(g,t):=(g-g_0)^4+(t-t_0)^4,
\quad
\chi_\varepsilon:=\chi+\varepsilon P.
\]
Since $v\le u^*\le\chi$ for every $v\in\mathcal S$,
\[
v-\chi_\varepsilon\le-\varepsilon P
\longrightarrow-\infty\quad\text{as }|g|\to\infty,
\]
uniformly in $t$ and $v\in\mathcal S$.
Choose $v_n\in\mathcal S$ and $(\hat g_n,\hat t_n)\to(g_0,t_0)$ such that
\[
v_n(\hat g_n,\hat t_n)\longrightarrow u^*(g_0,t_0).
\]
Let $(g_n,t_n)$ maximize $v_n-\chi_\varepsilon$ on $Q$, and set
\[
\gamma_n:=\max_Q(v_n-\chi_\varepsilon).
\]
Then
\[
0\ge \gamma_n\ge v_n(\hat g_n,\hat t_n)
-\chi_\varepsilon(\hat g_n,\hat t_n)\longrightarrow0,
\quad
0\le\varepsilon P(g_n,t_n)\le-\gamma_n\longrightarrow0.
\]
Therefore,
\begin{equation}\label{eq:perron_contact_sequence}
(g_n,t_n)\to(g_0,t_0),\quad
v_n(g_n,t_n)=\chi_\varepsilon(g_n,t_n)+\gamma_n
\longrightarrow u^*(g_0,t_0).
\end{equation}
For all sufficiently large $n$, $t_n<T$. The normalized upper test is
$\chi_\varepsilon+\gamma_n.$
Thus
\[
\min\{F_{\chi_\varepsilon}(g_n,t_n;v_n),B[v_n](g_n,t_n)\}\le0.
\]

For $\zeta\in\mathcal T$, the Marchaud operator is continuous at each
point with $t<T$. Indeed, on a compact set $K\subset\R\times[0,T)$,
for sufficiently small $d>0$,
\[
\sup_{(g,t)\in K}\int_0^d
\frac{|\zeta(g,t)-\zeta(g,t+\tau)|}{\tau^{1+\beta}}\,d\tau
\le\frac{C_Kd^{1-\beta}}{1-\beta}\longrightarrow0
\quad(0<\beta<1).
\]
The terminal term and the integral away from $\tau=0$ are continuous by
continuity of $\zeta$ on the closed strip. For $\beta=1$, continuity
follows from that of $\zeta_t$.
Using this fact and
\[
\limsup_{n\to\infty}v_n(-g_n,t_n)\le u^*(-g_0,t_0),
\]
we obtain, for fixed $\varepsilon$,
\begin{align*}
F_{\chi_\varepsilon}(g_n,t_n;v_n)
&\longrightarrow F_{\chi_\varepsilon}(g_0,t_0;u^*),\quad\text{and}\quad
\liminf_{n\to\infty}B[v_n](g_n,t_n)\ge B[u^*](g_0,t_0).
\end{align*}
Hence
\[
\min\{F_{\chi_\varepsilon}(g_0,t_0;u^*),B[u^*](g_0,t_0)\}
\le\liminf_{n\to\infty}
\min\{F_{\chi_\varepsilon}(g_n,t_n;v_n),B[v_n](g_n,t_n)\}\le0.
\]
Finally, $P=P_g=P_{gg}=0$ at $(g_0,t_0)$ and
$\rCder{\beta}_tP(g_0,t_0)$ is finite, so
\[
F_{\chi_\varepsilon}(g_0,t_0;u^*)
=F_\chi(g_0,t_0;u^*)
+\varepsilon\rCder{\beta}_tP(g_0,t_0)
\longrightarrow F_\chi(g_0,t_0;u^*).
\]
Letting $\varepsilon\to0$ therefore gives
\[
\min\{F_\chi(g_0,t_0;u^*),B[u^*](g_0,t_0)\}\le0.
\]
Together with \eqref{eq:perron_terminal_bound}, this proves that $u^*$
is a subsolution with terminal value zero. Thus $u^*\in\mathcal S$, and
\[
u\le u^*\le\sup_{w\in\mathcal S}w=u.
\]
Therefore $u=u^*\in\mathcal S$.

Before treating a failure of the differential branch, we verify the
pointwise obstacle condition required of a supersolution.  Suppose that
$B[u_*](p_0)<0$ at some $p_0=(g_0,t_0)\in Q_0$.  Let $p_n$ minimize
$u_*(p)+n|p-p_0|^2$ over $Q$, and define
\[
q_n(p):=u_*(p_n)+n|p_n-p_0|^2-n|p-p_0|^2.
\]
Then $q_n$ is a global lower test at $p_n$, $p_n\to p_0$, and
$u_*(p_n)\to u_*(p_0)$.  Lower semicontinuity at the reflected point gives
\[
\limsup_{n\to\infty}B[u_*](p_n)\le B[u_*](p_0)<0.
\]
For large $n$, the test $q_n$ therefore yields a strict failure of the
supersolution inequality.  Applying the bump argument below to this test
produces a member of $\mathcal S$ strictly above $u$, a contradiction.
Thus $B[u_*]\ge0$ on $Q_0$.

We argue by contradiction. Suppose there exist $(g_0,t_0)\in Q$, $t_0<T$, a global lower test $\varphi_0\in\mathcal T$, and $\theta>0$ such that
\[
u_*-\varphi_0\ge0\text{ on }Q,\quad
u_*(g_0,t_0)=\varphi_0(g_0,t_0),
\]
and
\begin{equation}\label{eq:perron_fail}
\min\{F_{\varphi_0}(g_0,t_0;u_*),B[u_*](g_0,t_0)\}
\le-4\theta.
\end{equation}
To make the contact strict and localize the lift on the unbounded strip,
set
\[
p(g,t):=(g-g_0)^4+(t-t_0)^4,
\quad
\varphi:=\varphi_0-\varepsilon p.
\]
At $(g_0,t_0)$, the perturbation and its first two spatial derivatives
vanish. Moreover,
\[
F_\varphi(g_0,t_0;u_*)
=F_{\varphi_0}(g_0,t_0;u_*)+\varepsilon K_\beta,
\quad
K_\beta:=
\begin{cases}
\displaystyle\frac{4(T-t_0)^{4-\beta}}
{(4-\beta)\Gamma(1-\beta)},&0<\beta<1,\\[1ex]
0,&\beta=1.
\end{cases}
\]
Choose $\varepsilon>0$ so that $\varepsilon K_\beta\le\theta$.
Then
\begin{equation}\label{eq:perron_strict_failure}
\min\{F_\varphi(g_0,t_0;u_*),B[u_*](g_0,t_0)\}
\le-3\theta,
\quad
u_*-\varphi\ge\varepsilon p\text{ on }Q.
\end{equation}

First, the contact point lies strictly below the upper barrier. Indeed,
if $u_*(g_0,t_0)=\tilde u(g_0,t_0)$, then
$\varphi\le u_*\le\tilde u$ globally, with equality at $(g_0,t_0)$.
The supersolution property of $\tilde u$ gives
\[
F_\varphi(g_0,t_0;u_*)
=F_\varphi(g_0,t_0;\tilde u)\ge0,
\]
whereas
\[
B[u_*](g_0,t_0)
\ge\tilde u(g_0,t_0)-\tilde u(-g_0,t_0)-\psi(g_0,t_0)
=c>0.
\]
This contradicts \eqref{eq:perron_strict_failure}. Hence
\[
m_1:=\tilde u(g_0,t_0)-u_*(g_0,t_0)>0.
\]

At least one of the two branches in \eqref{eq:perron_strict_failure}
is at most $-3\theta$. The function
$(g,t)\mapsto F_\varphi(g,t;\varphi)$ is continuous near $(g_0,t_0)$,
and
$(g,t)\mapsto\varphi(g,t)-u_*(-g,t)-\psi(g,t)$ is upper semicontinuous.
Thus, for sufficiently small $\delta\in(0,T-t_0)$, on
\[
Q_\delta:=\{(g,t)\in Q:(g-g_0)^2+(t-t_0)^2<\delta^2\},
\]
we have
\begin{equation}\label{eq:perron_local_gap}
\tilde u-\varphi\ge\frac{m_1}{2},
\end{equation}
and, according to the failing branch, either
\begin{align}
F_\varphi(g,t;\varphi)&\le-2\theta
&&\text{throughout }Q_\delta,\label{eq:perron_case_F}\\
\intertext{or}
\varphi(g,t)-u_*(-g,t)-\psi(g,t)&\le-2\theta
&&\text{throughout }Q_\delta.\label{eq:perron_case_B}
\end{align}
Since $a^4+b^4\ge\tfrac12(a^2+b^2)^2$, the global estimate in
\eqref{eq:perron_strict_failure} gives
\begin{equation}\label{eq:perron_exterior_gap}
u_*-\varphi\ge m_2:=\frac{\varepsilon\delta^4}{2}>0
\quad\text{on }Q\setminus Q_\delta.
\end{equation}
Choose
\[
0<\mu<\min\left\{\frac{m_1}{2},m_2,\theta,\frac{\theta}{r}\right\},
\quad
\Phi:=\varphi+\mu,
\quad
w:=\max\{u,\Phi\}.
\]
By \eqref{eq:perron_local_gap}--\eqref{eq:perron_exterior_gap},
\[
w=u\text{ on }Q\setminus Q_\delta,
\quad
\underline u\le u\le w\le\tilde u\text{ on }Q.
\]
In particular, $w\in\operatorname{USC}(Q)$ and $w(g,T)=0$.

Let $\xi\in\mathcal T$ be a normalized global upper test for $w$ at
$(g_1,t_1)$, with $t_1<T$.

\emph{If $w(g_1,t_1)=u(g_1,t_1)$}, then $\xi\ge w\ge u$ on $Q$
with equality at $(g_1,t_1)$. Thus $\xi$ is also an upper test for $u$.
Moreover,
\[
F_\xi(g_1,t_1;w)=F_\xi(g_1,t_1;u),\quad
B[w](g_1,t_1)\le B[u](g_1,t_1).
\]
The subsolution property of $u$ therefore yields
\[
\min\{F_\xi(g_1,t_1;w),B[w](g_1,t_1)\}\le0.
\]

\emph{If $w(g_1,t_1)=\Phi(g_1,t_1)>u(g_1,t_1)$}, then
$(g_1,t_1)\in Q_\delta$ and
\[
\xi\ge w\ge\Phi\text{ on }Q,
\quad
\xi(g_1,t_1)=\Phi(g_1,t_1).
\]
Therefore,
\[
\xi_g(g_1,t_1)=\Phi_g(g_1,t_1),\quad
\xi_{gg}(g_1,t_1)\ge\Phi_{gg}(g_1,t_1).
\]
For $0<\beta<1$, writing $h(s):=\xi(g_1,s)-\Phi(g_1,s)\ge0$
with $h(t_1)=0$, the Marchaud formula gives
\begin{align*}
\rCder{\beta}_t\xi(g_1,t_1)-\rCder{\beta}_t\Phi(g_1,t_1)
&=-\frac{h(T)}{\Gamma(1-\beta)(T-t_1)^\beta}
-\frac{\beta}{\Gamma(1-\beta)}
\int_{t_1}^T\frac{h(s)}{(s-t_1)^{1+\beta}}\,ds\\
&\le0.
\end{align*}
For $\beta=1$, the corresponding inequality follows from
$-\xi_t(g_1,t_1)\le-\Phi_t(g_1,t_1)$, with the right derivative at
$t_1=0$. Hence
\begin{align*}
F_\xi(g_1,t_1;w)
&\le F_\Phi(g_1,t_1;\Phi)
 =F_\varphi(g_1,t_1;\varphi)+r\mu,\\
B[w](g_1,t_1)
&\le\varphi(g_1,t_1)-u_*(-g_1,t_1)-\psi(g_1,t_1)+\mu.
\end{align*}
If \eqref{eq:perron_case_F} holds, then
\[
F_\xi(g_1,t_1;w)\le-2\theta+r\mu<-\theta.
\]
If \eqref{eq:perron_case_B} holds, then
\[
B[w](g_1,t_1)\le-2\theta+\mu<-\theta.
\]
In either case the required minimum is nonpositive. Thus $w\in\mathcal S$.

By the definition of the Perron envelope, $u\ge w\ge\Phi$ on $Q$.
Taking the lower limit at $(g_0,t_0)$ gives
\[
u_*(g_0,t_0)
\ge\Phi(g_0,t_0)
=\varphi(g_0,t_0)+\mu
=u_*(g_0,t_0)+\mu,
\]
a contradiction. Therefore $u_*$ is a viscosity supersolution.

By \eqref{eq:perron_terminal_bound}--\eqref{eq:perron_envelope_bounds},
\[
u^*(g,T)=u_*(g,T)=0.
\]
The semicontinuous comparison theorem applies to $u=u^*$ and $u_*$,
yielding
\[
u=u^*\le u_*\le u.
\]
Hence $u=u^*=u_*$ is continuous and solves \eqref{eq:pde}, with
$\underline u\le u\le\tilde u$. Applying the same comparison principle
in both directions to any two continuous solutions of at most linear growth
proves uniqueness.
\end{proof}

\begin{corollary}
\label{cor:contact-maturity}
Let $u$ be the solution of Theorem~\ref{thm:exist}, and let $f$
be the profile in Lemma~\ref{lem:super}. Define
\[
\mathcal E:=\{(g,t)\in\R\times[0,T):B[u](g,t)=0\}.
\]
Then
\[
\mathcal E\subseteq
\{(g,t)\in\R\times[0,T):\alpha(t)f(g)\ge c\}.
\]
\end{corollary}

\begin{proof}
The barrier bounds and the definition of $f$ imply
\[
0\le u(g,t)\le\alpha(t)f(g),
\quad
f(g)-f(-g)=g.
\]
Since $\psi(g,t)=g\alpha(t)-c$, we obtain
\begin{align*}
B[u](g,t)
&=u(g,t)-u(-g,t)-g\alpha(t)+c
\ge-\alpha(t)f(-g)-g\alpha(t)+c\\
&=c-\alpha(t)\bigl(f(-g)+g\bigr)=c-\alpha(t)f(g).
\end{align*}
Thus $B[u](g,t)=0$ implies $\alpha(t)f(g)\ge c$,
proving the asserted inclusion.
\end{proof}

\begin{proposition}
Under the hypotheses of Theorem~\ref{thm:exist}, let $u_c$ denote the
solution with transaction cost $c>0$, keeping all other data
fixed. Then $u_c$ is non-increasing in $c$, and \[
u_c\longrightarrow 0
\quad\text{locally uniformly on }\mathbb R\times[0,T]
\quad\text{as }c\to\infty.
\]
On each bounded spatial interval the convergence has rate $O(1/c)$,
uniformly with respect to $\beta\in(0,1]$.
\end{proposition}
\begin{proof}
Let $c_1\le c_2.$ Writing
\[
B_c[w](g,t):=w(g,t)-w(-g,t)-g\alpha(t)+c,
\]
we have
\[
B_{c_2}[u_{c_1}]
=B_{c_1}[u_{c_1}]+c_2-c_1\ge0.
\]
The differential supersolution inequality is unchanged.
Thus $u_{c_1}$ is a supersolution of the problem with cost $c_2$.
Since both solutions vanish at maturity,
Theorem~\ref{thm:comp} yields $u_{c_2}\le u_{c_1}$.
Hence increasing transaction costs decreases the speculative
premium.

If $\alpha(0)=0$, then $\alpha\equiv0$ and $u_c\equiv0$.
Otherwise, set $b=2c/\alpha(0)$, $K=5\sigma^2/(4r)$, and
$$
U(g,t)=\alpha(t)\left(
\frac{g^2}{2b}+\frac{g^3}{2(b^2+g^2)}+\frac Kb\right).
$$
Using $ b|g|/(b^2+g^2)\le 1/2,$ we obtain
\begin{equation}\label{eq:cost-barrier-bounds}
\alpha(t)\left(\frac{g^2}{4b}+\frac Kb\right)
\le U(g,t)\le
\alpha(t)\left(\frac{3g^2}{4b}+\frac Kb\right).
\end{equation}
In particular, $U\ge0$. Direct differentiation gives
\begin{align*}
gU_g(g,t)
&=\frac{\alpha(t)g^2}{b}
\left(1+\frac{bg(3b^2+g^2)}{2(b^2+g^2)^2}\right)
\ge\frac{\alpha(t)g^2}{4b},\\
U_{gg}(g,t)
&=\alpha(t)\left(
\frac1b+\frac{b^2g(3b^2-g^2)}{(b^2+g^2)^3}
\right)
\le\frac{5\alpha(t)}{2b}.
\end{align*}
Here we used
\[
\left|\frac{bg(3b^2+g^2)}{2(b^2+g^2)^2}\right|
\le\frac34,
\quad
\left|\frac{b^2g(3b^2-g^2)}{(b^2+g^2)^3}\right|
\le\frac{3}{2b}.
\]
Since $\alpha'\le0$ and the spatial factor in $U$ is nonnegative,
$\rCder{\beta}_t U\ge0$ for every $\beta\in(0,1]$.
Therefore,
\[
F[U]\ge
\frac{\alpha(t)}b
\left(\frac{\rho+r}{4}g^2+rK-\frac{5\sigma^2}{4}\right)
\ge0.
\]
On the other hand,
$$
B_c[U]=c-\alpha(t)\frac{gb^2}{b^2+g^2}
\ge c-\frac{\alpha(t)b}{2}\ge0.
$$
Since $U\ge 0$ and $U(\cdot,T)=0$, Theorem~\ref{thm:comp} gives
$0\le u_c\le U$. Hence, for every $R>0$,
$$
\|u_c\|_{L^\infty([-R,R]\times[0,T])}
\le\frac{\alpha(0)^2}{2c}
   \left(\frac{3R^2}{4}+K\right)
=O_R(c^{-1}),
$$
with a constant independent of $\beta$.
\end{proof}

\section{Regularity of Solutions}\label{sec:regularity}

In this section we establish spatial Lipschitz continuity and $\beta/2$-regularity in time  of solutions to \eqref{eq:pde}.
To prove Lipschitz regularity, we first center the obstacle and show that the symmetric envelope $V_a(g,t) := \max_{|h|\le a} v(g+h,t) - La$ (see \eqref{eq:sym} for the definition of $v$) is a subsolution for sufficiently large $L$.
The comparison principle yields $V_a\le v$; applying this
inequality to opposite shifts and restoring the linear profile gives
the desired Lipschitz estimate for $u$.  This idea is carried out in the proof of the following proposition.

\begin{proposition}[Spatial Lipschitz continuity]
\label{prop:spatial-lipschitz}
Under Assumption~\ref{ass:standing}, suppose that
$\alpha\in C^1([0,T]),$ also $ \alpha'(t)\le 0$ and $\alpha(T)=0.$
Let $u\in C(\R\times[0,T])$ be the viscosity solution of
\eqref{eq:pde}, with uniform linear growth and $u(g,T)=0$.
Then there exists $C<\infty$, depending only on the data, such that
\[
 |u(g+h,t)-u(g,t)|\le C|h|,
 \quad g,h\in\R,\quad t\in[0,T].
\]
\end{proposition}

\begin{proof}
Define the smooth profile $\psi_0(g,t):=\frac{g}{2}\alpha(t)$ and the centered unknown
\begin{equation}\label{eq:sym}
     v(g,t):=u(g,t)-\psi_0(g,t)=u(g,t)-\frac{g}{2}\alpha(t),\quad
 l(t):=\frac12\bigl(\rCder{\beta}_t\alpha+(r+\rho)\alpha\bigr)(t).
\end{equation}
Since $\alpha\in C^1([0,T])$, its right fractional derivative $\rCder{\beta}_t\alpha$ is continuous on $[0,T]$ for all $\beta\in(0,1].$ Moreover, it is bounded by 
$$
\|\alpha'\|_\infty \frac{T^{1-\beta}}{\Gamma(2-\beta)},\quad \text{for}\quad 0<\beta<1,
$$
and equal to $-\alpha'$ for $\beta=1$. Thus $l\in C([0,T])$ and $\|l\|_\infty<\infty$.

Let us verify that $v$ is a viscosity solution, in the sense of Definition~\ref{def:visc}, of the following obstacle problem
\begin{equation}\label{eq:centered-pde}
 \min\Bigl\{F[v](g,t)+g\,l(t),\;\; v(g,t)-v(-g,t)+c\Bigr\}=0,
 \quad v(g,T)=0.
\end{equation}
First, at maturity $t=T$, $\alpha(T)=0$ and $u(g,T)=0$ imply $v(g,T)=0$. Second, because $\psi_0(-g,t)=-\psi_0(g,t)$ and $\psi(g,t)=g\alpha(t)-c$, we obtain
\[
 u(g,t)-u(-g,t)-\psi(g,t)
 =v(g,t)-v(-g,t)+c.
\]
 Third, $\varphi$ is a normalized local upper (resp.\ lower) test function for $v$ at $(g_0,t_0)$ with $t_0<T$ if and only if $\tilde\varphi:=\varphi+\psi_0$ is a normalized local upper (resp.\ lower) test function for $u=v+\psi_0$ at $(g_0,t_0)$. Because $\psi_0\in C(\R\times[0,T])\cap C^{2,1}(\R\times[0,T))$ is smooth and linear in $g$, its near-time test integral on $(0,\delta)$, distant-time integral on $(\delta,T-t_0)$, and terminal term in the split operator $D_\delta$ combine into the exact Marchaud representation of $\rCder{\beta}_t\psi_0(g_0,t_0)=\frac{g_0}{2}\rCder{\beta}_t\alpha(t_0)$ for every admissible $\delta>0$. Therefore,
\begin{equation}\label{eq: equiv-visc}
\begin{gathered}
 F_{\tilde\varphi,\delta}(g_0,t_0;u)
 =D_\delta[v,\varphi](g_0,t_0)+\frac{g_0}{2}\rCder{\beta}_t\alpha(t_0)
   +\rho g_0\Bigl(\varphi_g(g_0,t_0)+\frac12\alpha(t_0)\Bigr)\\
   -\frac{\sigma^2}{2}\varphi_{gg}(g_0,t_0)
   +r\Bigl(v(g_0,t_0)+\frac{g_0}{2}\alpha(t_0)\Bigr)
 =F_{\varphi,\delta}(g_0,t_0;v)+g_0\,l(t_0).
\end{gathered}
\end{equation}

The identical relation holds for $\beta=1$. Thus $u$ is a viscosity subsolution (resp.\ supersolution) of \eqref{eq:pde} if and only if $v$ is a viscosity subsolution (resp.\ supersolution) of \eqref{eq:centered-pde}.

By Assumption~\ref{ass:standing}, $r>0$ and $\rho\ge0$. Fix a constant
\[
 L\ge\frac{\|l\|_\infty}{r},
\]
and, for each $a>0$, define on $Q=\R\times[0,T]$
\[
 V_a(g,t):=\max_{|h|\le a}v(g+h,t)-La.
\]
The maximum is attained because $v$ is continuous and $[-a,a]$ is compact. The same compactness gives $V_a\in C(Q)$. Moreover, if $|v(g,t)|\le C_1(1+|g|)$, then
\[
 |V_a(g,t)|\le C_1(1+|g|+a)+La,
\]
so $V_a$ has uniform linear growth for each fixed $a>0$. At maturity,
\[
 V_a(g,T)=-La\le0.
\]

We claim that $V_a$ is a viscosity subsolution of the  problem \eqref{eq:centered-pde}. Let $(g_0,t_0)\in\R\times[0,T)$ and let $\varphi$ be a $C^{2,1}$ test function touching $V_a$ locally from above at $(g_0,t_0)$ in a neighborhood $U\subset\R\times[0,T)$, normalized so that
\[
 \varphi(g_0,t_0)=V_a(g_0,t_0),\quad \varphi(g,s)\ge V_a(g,s)\quad\text{for all }(g,s)\in U.
\]
Let $h_0\in[-a,a]$ be a shift attaining the maximum in $V_a(g_0,t_0)$, and set
\[
 y_0:=g_0+h_0,\quad d:=La\ge0,
 \quad\text{so that}\quad
 V_a(g_0,t_0)=v(y_0,t_0)-d.
\]
Define the spatially translated test function
\[
 \chi(z,s):=\varphi(z-h_0,s)+d
\]
on the shifted neighborhood $U+(h_0,0)$ of $(y_0,t_0)$. For every $(z,s)\in U+(h_0,0)$, the definition of $V_a(z-h_0,s)$ with the candidate shift $h_0$ gives
\[
 v(z,s)=v\bigl((z-h_0)+h_0,s\bigr)
 \le V_a(z-h_0,s)+d
 \le \varphi(z-h_0,s)+d
 =\chi(z,s),
\]
with equality at $(z,s)=(y_0,t_0)$ since $\chi(y_0,t_0)=\varphi(g_0,t_0)+d=V_a(g_0,t_0)+d=v(y_0,t_0)$. Thus $\chi$ is a normalized local upper test function for $v$ at $(y_0,t_0)$, and its derivatives satisfy
\begin{equation}\label{eq:chi-derivs}
 \chi_g(y_0,t_0)=\varphi_g(g_0,t_0),\quad
 \chi_{gg}(y_0,t_0)=\varphi_{gg}(g_0,t_0),\quad
 \chi_t(y_0,t_0)=\varphi_t(g_0,t_0).
\end{equation}

Furthermore, for $0<s<1$ sufficiently small, the shift $(1-s)h_0$ is admissible at $g_0+s h_0$, since
\[
 |(1-s)h_0|\le a,\quad
 (g_0+s h_0)+(1-s)h_0=y_0.
\]
Keeping $y_0$ fixed therefore gives
\[
 \varphi(g_0+s h_0,t_0)
 \ge V_a(g_0+s h_0,t_0)
 \ge v(y_0,t_0)-d
 =\varphi(g_0,t_0).
\]
Subtracting $\varphi(g_0,t_0)$, dividing by $s>0$, and letting $s\downarrow0$, we obtain
\begin{equation}\label{eq:opt-shift}
 h_0\varphi_g(g_0,t_0)\ge0.
\end{equation}

At the reflected spatial coordinate $-g_0$, testing the maximum in $V_a(-g_0,t_0)$ with the reflected shift $-h_0$ yields
\[
 V_a(-g_0,t_0)
 \ge v(-g_0-h_0,t_0)-La
 =v(-y_0,t_0)-d.
\]
Combining this with $V_a(g_0,t_0)=v(y_0,t_0)-d$ cancels $d$ and gives the obstacle comparison
\begin{equation}\label{eq:obs-shift}
 V_a(g_0,t_0)-V_a(-g_0,t_0)+c
 \le v(y_0,t_0)-v(-y_0,t_0)+c.
\end{equation}

Next, for any future time $s\in[t_0,T]$, testing the maximum in $V_a(g_0,s)$ with the fixed shift $h_0$ yields
\[
 V_a(g_0,s)\ge v(g_0+h_0,s)-d=v(y_0,s)-d.
\]
Subtracting this inequality from $V_a(g_0,t_0)=v(y_0,t_0)-d$ shows that for all $s\in[t_0,T]$,
\[
 V_a(g_0,t_0)-V_a(g_0,s)\le v(y_0,t_0)-v(y_0,s).
\]
At the same time, for every admissible $\delta\in(0,T-t_0)$ and every $\tau\in[0,\delta]$, the definition of $\chi$ gives
\[
 \varphi(g_0,t_0)-\varphi(g_0,t_0+\tau)=\chi(y_0,t_0)-\chi(y_0,t_0+\tau).
\]
Since the kernel $\tau^{-1-\beta}$ and the terminal weight $(T-t_0)^{-\beta}$ in the split operator $D_\delta$ are nonnegative, integrating these inequalities over $\tau\in(\delta,T-t_0)$ and $\tau\in(0,\delta)$ yields, for $0<\beta<1$,
\begin{equation}\label{eq:hist-shift}
 D_\delta[V_a,\varphi](g_0,t_0)\le D_\delta[v,\chi](y_0,t_0).
\end{equation}
For $\beta=1$, \eqref{eq:chi-derivs} gives $-\varphi_t(g_0,t_0)=-\chi_t(y_0,t_0)$, so \eqref{eq:hist-shift} holds with equality.

We now verify that at $(g_0,t_0)$, for every admissible $\delta$,
\begin{equation}\label{eq:Va-sub}
 \min\Bigl\{F_{\varphi,\delta}(g_0,t_0;V_a)+g_0\,l(t_0),\;\;
 V_a(g_0,t_0)-V_a(-g_0,t_0)+c\Bigr\}\le0.
\end{equation}
If $V_a(g_0,t_0)-V_a(-g_0,t_0)+c\le0$, then \eqref{eq:Va-sub} holds immediately. Suppose instead that
\[
 V_a(g_0,t_0)-V_a(-g_0,t_0)+c>0.
\]
Then \eqref{eq:obs-shift} implies $v(y_0,t_0)-v(-y_0,t_0)+c>0$. Because $v$ is a viscosity subsolution of \eqref{eq:centered-pde} and $\chi$ touches $v$ locally from above at $(y_0,t_0)$, the active subsolution inequality for $v$ at $(y_0,t_0)$ implies
\begin{equation}\label{eq:v-diff-active}
 F_{\chi,\delta}(y_0,t_0;v)+y_0\,l(t_0)\le0.
\end{equation}
Using \eqref{eq:chi-derivs}, \eqref{eq:hist-shift}, $g_0=y_0-h_0$, and $V_a(g_0,t_0)=v(y_0,t_0)-La$, we obtain
\begin{align*}
 F_{\varphi,\delta}(g_0,t_0;V_a)+g_0\,l(t_0)
 &\le F_{\chi,\delta}(y_0,t_0;v)+y_0\,l(t_0)
 -\rho h_0\varphi_g(g_0,t_0)-rLa-h_0\,l(t_0)\\
 &\le-rLa+|h_0|\,\|l\|_\infty
 \le a\bigl(\|l\|_\infty-rL\bigr)\le0.
\end{align*}
Here we used \eqref{eq:v-diff-active}, \eqref{eq:opt-shift}, $\rho\ge0$, and $|h_0|\le a$. Thus \eqref{eq:Va-sub} holds in both cases, proving that $V_a$ is a viscosity subsolution of \eqref{eq:centered-pde}.

By \eqref{eq: equiv-visc}, the function $U_a(g,t):=V_a(g,t)+\frac{g}{2}\alpha(t)$ is a viscosity subsolution of the original obstacle problem \eqref{eq:pde} in the uniform linear-growth class, with $U_a(g,T)=-La\le0$. Applying the Comparison Principle (Theorem~\ref{thm:comp}) to $U_a$ and $u$ yields $U_a\le u$ on $\R\times[0,T]$, or equivalently $V_a\le v$ on $\R\times[0,T]$.

Therefore, for every $a>0$, $|h|\le a$, $g\in\R$, and $t\in[0,T]$,
\[
 v(g+h,t)-La\le V_a(g,t)\le v(g,t).
\]
For $h\ne0$, take $a=|h|$ to obtain $v(g+h,t)-v(g,t)\le L|h|$; the same inequality is immediate for $h=0$. Rewriting this in terms of $u(g,t)=v(g,t)+\frac{g}{2}\alpha(t)$, we get
\[
 u(g+h,t)-u(g,t)
 =v(g+h,t)-v(g,t)+\frac{h}{2}\alpha(t)
 \le\left(L+\frac12\|\alpha\|_\infty\right)|h|.
\]
Set $C:=L+\frac12\|\alpha\|_\infty$. We have proved that
\[
 u(g+h,t)-u(g,t)\le C|h|
 \quad\text{for all }g,h\in\R,\quad t\in[0,T].
\]
Applying this inequality with base point $g+h$ and increment $-h$
gives
\[
 u((g+h)-h,t)-u(g+h,t)\le C|-h|,
\]
or equivalently,
\[
 u(g,t)-u(g+h,t)\le C|h|.
\]
Thus,
\[
 -C|h|\le u(g+h,t)-u(g,t)\le C|h|,
\]
and therefore
\[
 |u(g+h,t)-u(g,t)|\le C|h|,
 \quad g,h\in\R,\quad t\in[0,T].
\]
This proves the asserted spatial Lipschitz continuity.
\end{proof}

Next, we  establish temporal H\"older continuity using the spatial
Lipschitz estimate. The reflected obstacle structure gives local bounds
for the fractional heat operator, which we transfer to the distributional
setting. Spatial mollification then allows us to estimate time increments
through the fractional-integral representation. Balancing this estimate
with the mollification error yields temporal $C^{\beta/2}$ continuity,
uniformly on bounded spatial intervals and throughout $[0,T]$.

\begin{proposition}[Temporal H\"older continuity]
\label{prop:temporal-holder}
Under the assumptions of Proposition~\ref{prop:spatial-lipschitz},
let $u$ be the continuous viscosity solution of \eqref{eq:pde},
with $u(g,T)=0$ and uniform linear growth.
Then, for every $R>0$, there exists $C_R<\infty$ such that
\begin{equation}\label{eq:temporal-holder}
 |u(g,t)-u(g,s)|\le C_R|t-s|^{\beta/2},
 \quad |g|\le R,\quad s,t\in[0,T].
\end{equation}
The constant depends only on $R,T,\beta,r,\rho,\sigma,c$,
$\|\alpha\|_{C^1([0,T])}$, and the spatial Lipschitz constant
provided by Proposition~\ref{prop:spatial-lipschitz}.
\end{proposition}

\begin{proof}
Let $L<\infty$ be the global spatial Lipschitz constant of $u$ provided by Proposition~\ref{prop:spatial-lipschitz}.
We first show that on $\R\times[0,T)$,
\begin{equation}\label{eq:time-linear-bounds}
 0\le F[u]\le \mathcal H_+
 \quad\hbox{in the viscosity sense},
\end{equation}
where
\[
 \mathcal H(g,t):=F[\psi](g,t)
 =g\bigl(\rCder{\beta}_t\alpha+(r+\rho)\alpha\bigr)(t)-rc,
\]
and $\mathcal H_+(g,t) := \max\{\mathcal H(g,t),\, 0\}$.

Since $u$ is a viscosity supersolution of \eqref{eq:pde}, the lower inequality $F_{\varphi,\delta}(g_0,t_0;u)\ge0$ holds at every normalized local lower contact $(g_0,t_0)$ of $u-\varphi$, establishing $0\le F[u]$.

To prove the upper inequality $F[u]\le \mathcal H_+$, let $(g_0,t_0)\in\R\times[0,T)$ and let $\varphi$ be a $C^{2,1}$ test function touching $u$ locally from above at $(g_0,t_0)$ on a neighborhood $U\subset\R\times[0,T)$, normalized so that $\varphi(g_0,t_0)=u(g_0,t_0)$ and $\varphi\ge u$ on $U$. If the obstacle is strictly inactive at $(g_0,t_0)$, that is, $B[u](g_0,t_0)=u(g_0,t_0)-u(-g_0,t_0)-\psi(g_0,t_0)>0$, then the viscosity subsolution inequality for $u$ directly gives
\[
 F_{\varphi,\delta}(g_0,t_0;u)\le 0\le \mathcal H_+(g_0,t_0)
\]
for every admissible $\delta$. Suppose instead that $B[u](g_0,t_0)\le 0$. Because $u$ is also a viscosity supersolution, it satisfies the pointwise obstacle constraint $B[u]\ge0$ everywhere on $\R\times[0,T)$. Hence $B[u](g_0,t_0)=0$, which means
\begin{equation}\label{eq:contact-obs-eq}
 u(g_0,t_0)=u(-g_0,t_0)+\psi(g_0,t_0),
\end{equation}
and the complementarity identity $B[u](g_0,t_0)+B[u](-g_0,t_0)=2c$ yields
\[
 B[u](-g_0,t_0)=2c>0.
\]
Define the reflected smooth function
\[
 \chi(x,s):=\varphi(-x,s)-\psi(-x,s)
\]
on the reflected neighborhood $\widetilde U:=\{(-x,s):(x,s)\in U\}$ of $(-g_0,t_0)$. For every $(x,s)\in\widetilde U$, combining the global obstacle inequality $B[u](-x,s)\ge0$ with the upper contact inequality $u(-x,s)\le\varphi(-x,s)$ on $U$ gives
\[
 u(x,s)\le u(-x,s)-\psi(-x,s)\le \varphi(-x,s)-\psi(-x,s)=\chi(x,s),
\]
while at $(x,s)=(-g_0,t_0)$, \eqref{eq:contact-obs-eq} gives
\[
 \chi(-g_0,t_0)=\varphi(g_0,t_0)-\psi(g_0,t_0)=u(g_0,t_0)-\psi(g_0,t_0)=u(-g_0,t_0).
\]
Thus $\chi$ is a normalized local upper test function for $u$ at $(-g_0,t_0)$ on $\widetilde U$, and any $\delta\in(0,T-t_0)$ admissible for $\varphi$ at $(g_0,t_0)$ is admissible for $\chi$ at $(-g_0,t_0)$. Moreover, since $B[u](-g_0,t_0)=2c>0$, the subsolution condition for $u$ at $(-g_0,t_0)$ implies 
\begin{equation}\label{eq:chi-sub-active}
 F_{\chi,\delta}(-g_0,t_0;u)\le 0.
\end{equation}

We now compare $F_{\varphi,\delta}(g_0,t_0;u)$ with $F_{\chi,\delta}(-g_0,t_0;u)$. For all $s\in[t_0,T]$, the inequality $u(g_0,s)\ge u(-g_0,s)+\psi(g_0,s)$ holds (including at $s=T$, where $u(g_0,T)=0>-c=u(-g_0,T)+\psi(g_0,T)$). Subtracting this from \eqref{eq:contact-obs-eq} yields, for all $s\in[t_0,T]$,
\[
 u(g_0,t_0)-u(g_0,s)\le\bigl(u(-g_0,t_0)-u(-g_0,s)\bigr)+\bigl(\psi(g_0,t_0)-\psi(g_0,s)\bigr).
\]
Meanwhile, for near times $s=t_0+\tau$ with $\tau\in[0,\delta]$, the definition of $\chi$ gives the exact identity
\[
 \varphi(g_0,t_0)-\varphi(g_0,t_0+\tau)
 =\bigl(\chi(-g_0,t_0)-\chi(-g_0,t_0+\tau)\bigr)+\bigl(\psi(g_0,t_0)-\psi(g_0,t_0+\tau)\bigr).
\]
For $0<\beta<1$, integrating these relations against the positive Marchaud kernel $\tau^{-1-\beta}$ and terminal weight $(T-t_0)^{-\beta}$, and noting that the smooth profile $\psi(g_0,\cdot)$ satisfies
$$
D_\delta[\psi,\psi](g_0,t_0)=\rCder{\beta}_t\psi(g_0,t_0)=g_0\rCder{\beta}_t\alpha(t_0),
$$ 
we obtain
\begin{equation}\label{eq:Ddelta-reflected}
 D_\delta[u,\varphi](g_0,t_0)
 \le D_\delta[u,\chi](-g_0,t_0)+\rCder{\beta}_t\psi(g_0,t_0).
\end{equation}
For $\beta=1$, differentiating $\varphi(g_0,s)=\chi(-g_0,s)+\psi(g_0,s)$ at $s=t_0$ gives $$
-\varphi_t(g_0,t_0)=-\chi_t(-g_0,t_0)-\psi_t(g_0,t_0),
$$
so \eqref{eq:Ddelta-reflected} holds with equality. Furthermore, differentiating $\chi(x,s)=\varphi(-x,s)+x\alpha(s)+c$ in space at $x=-g_0$ gives
\[
 \chi_g(-g_0,t_0)=-\varphi_g(g_0,t_0)+\alpha(t_0),\quad
 \chi_{gg}(-g_0,t_0)=\varphi_{gg}(g_0,t_0),
\]
so that $\rho g_0\varphi_g(g_0,t_0)=\rho(-g_0)\chi_g(-g_0,t_0)+\rho g_0\alpha(t_0)$. Combining this with \eqref{eq:contact-obs-eq}, \eqref{eq:chi-sub-active}, and \eqref{eq:Ddelta-reflected} yields
\begin{align*}
 F_{\varphi,\delta}(g_0,t_0;u)
 &=D_\delta[u,\varphi](g_0,t_0)+\rho g_0\varphi_g(g_0,t_0)-\frac{\sigma^2}{2}\varphi_{gg}(g_0,t_0)+ru(g_0,t_0)\\
 &\le D_\delta[u,\chi](-g_0,t_0)+\rho(-g_0)\chi_g(-g_0,t_0)-\frac{\sigma^2}{2}\chi_{gg}(-g_0,t_0)+ru(-g_0,t_0)\\
 &\quad +\rCder{\beta}_t\psi(g_0,t_0)+\rho g_0\alpha(t_0)+r\psi(g_0,t_0)\\
 &=F_{\chi,\delta}(-g_0,t_0;u)+\mathcal H(g_0,t_0)\le \mathcal H(g_0,t_0)\le \mathcal H_+(g_0,t_0).
\end{align*}
This establishes \eqref{eq:time-linear-bounds}.

Fix $R>0$ and set $Q_R:=(-R-2,R+2)\times(0,T)$. Since $\alpha\in C^1([0,T])$, $\mathcal H_+$ is bounded on $Q_R$. The established barrier bounds $0\le u(g,t)\le\alpha(t)f(g)$, with $f$ as in Lemma~\ref{lem:super}, give
\[
 \|u\|_{L^\infty(Q_R)}
 \le \|\alpha\|_{L^\infty(0,T)}
       \max_{|g|\le R+2}f(g).
\]
Fixing the barrier parameters in terms of the equation's data, this bound has the dependence stated in the proposition.
Let $L$ be the spatial Lipschitz constant from
Proposition~\ref{prop:spatial-lipschitz}. If $\varphi$ touches $u$
from above at $(g_0,t_0)$, then, for all sufficiently small $h$, we get
\[
-L|h|
\le u(g_0+h,t_0)-u(g_0,t_0)
\le \varphi(g_0+h,t_0)-\varphi(g_0,t_0).
\]
Dividing by $h>0$ and letting $h\downarrow0$ gives
$\varphi_g(g_0,t_0)\ge-L$. Dividing by $h<0$ and letting
$h\uparrow0$ gives $\varphi_g(g_0,t_0)\le L$.
At a lower contact, the analogous inequalities
\[
\varphi(g_0+h,t_0)-\varphi(g_0,t_0)
\le u(g_0+h,t_0)-u(g_0,t_0)
\le L|h|
\]
give the same conclusion. Thus
\[
|\varphi_g(g_0,t_0)|\le L
\]
at every smooth upper or lower contact.

Fix $R>0$ and recall
$Q_R:=(-R-2,R+2)\times(0,T)$. Set
\[
M_R:=\sup_{[-R-2,R+2]\times[0,T]}|u|,
\quad
A_R:=\sup_{Q_R}\mathcal H_+.
\]
These constants are finite. In particular, the barrier bound gives
\[
M_R\le\|\alpha\|_{L^\infty(0,T)}
       \max_{|g|\le R+2}f(g).
\]
Hence, at every contact point in $Q_R$,
\[
\bigl|\rho g_0\varphi_g(g_0,t_0)+ru(g_0,t_0)\bigr|
\le \rho(R+2)L+rM_R.
\]
At an upper contact, \eqref{eq:time-linear-bounds} yields,
for every admissible $\delta$,
\begin{align*}
D_\delta[u,\varphi](g_0,t_0)
-\frac{\sigma^2}{2}\varphi_{gg}(g_0,t_0)
&\le \mathcal H_+(g_0,t_0)
   -\rho g_0\varphi_g(g_0,t_0)-ru(g_0,t_0)\\
&\le A_R+\rho(R+2)L+rM_R.
\end{align*}
At a lower contact, the other inequality in
\eqref{eq:time-linear-bounds} gives
\begin{align*}
D_\delta[u,\varphi](g_0,t_0)
-\frac{\sigma^2}{2}\varphi_{gg}(g_0,t_0)
&\ge-\rho g_0\varphi_g(g_0,t_0)-ru(g_0,t_0)\ge-\rho(R+2)L-rM_R.
\end{align*}
For $\beta=1$, the same estimates hold with
$D_\delta[u,\varphi]=-\varphi_t$.
Therefore, choosing
\[
C_R:=A_R+\rho(R+2)L+rM_R,
\]
we obtain
\begin{equation}\label{eq:time-heat-bounds}
-C_R\le\mathcal L_\beta u\le C_R
\quad\hbox{on }Q_R,
\quad
\mathcal L_\beta u
:=\rCder{\beta}_t u-\frac{\sigma^2}{2}u_{gg},
\end{equation}
in the viscosity sense.

To pass to distributions, extend $u$ by $0$ for $t>T$ and denote this extension by $\check u$. Let $U_n,V_n$ be its spatio-temporal sup- and inf-convolutions, with parameters tending to zero, as constructed in the proof of Lemma~\ref{lem:right-sums}. For $0<\beta<1$, use their full future histories
\[
 \mathcal D^\beta W(g,t)
 :=\frac{\beta}{\Gamma(1-\beta)}
   \int_0^\infty\frac{W(g,t)-W(g,t+\tau)}{\tau^{1+\beta}}\,d\tau,
\]
and set $\mathcal D^1W=-W_t$.

The histories of $U_n,V_n$ use their values also for $t>T$.
For $0<\beta<1$, set $c_\beta:=\beta/\Gamma(1-\beta)$.
Given $\xi\in C_c^\infty(Q_R)$, extend it by zero and define
\begin{align*}
 (\mathcal D^\beta)^*\xi(g,t)
 &:=c_\beta\int_0^\infty
       \frac{\xi(g,t)-\xi(g,t-\tau)}{\tau^{1+\beta}}\,d\tau,\\
 \langle\mathcal D^\beta W,\xi\rangle
 &:=\int_{\R}\int_0^\infty W(g,t)(\mathcal D^\beta)^*\xi(g,t)\,dt\,dg .
\end{align*}
This pairing is finite for $W$ bounded on each strip
$[-R',R']\times[0,\infty)$. Indeed, the adjoint has compact spatial
support and is bounded since
$|\xi(g,t)-\xi(g,t-\tau)|\le C_\xi\min\{\tau,1\}$. For $t>2T$ it satisfies
\[
 (\mathcal D^\beta)^*\xi(g,t)
 =-c_\beta\int_0^T\frac{\xi(g,s)}{(t-s)^{1+\beta}}\,ds
 =O_\xi(t^{-1-\beta})
\]
uniformly in $g$. For locally Lipschitz $W$ with these strip bounds,
truncating at $\tau>d$, changing variables, and letting $d\downarrow0$
shows that the pairing agrees with the full future history above.
Write $I_{0+}^{1-\beta}$ and $I_{T-}^{1-\beta}$ for the left and right
Riemann--Liouville integrals in time. Integration by parts in $\tau$ and
the identity $\xi(g,0)=0$ give
$(\mathcal D^\beta)^*\xi=I_{0+}^{1-\beta}\partial_t\xi
=\partial_t I_{0+}^{1-\beta}\xi$ for $t>0$.
Since $\check u=0$ for $t>T$, Fubini's theorem gives
\begin{align*}
 \langle\mathcal D^\beta\check u,\xi\rangle
 =\int_{Q_R}u\,I_{0+}^{1-\beta}\partial_t\xi\,dg\,dt
 =\int_{Q_R}(I_{T-}^{1-\beta}u)\,\partial_t\xi\,dg\,dt
 =\langle-\partial_t I_{T-}^{1-\beta}u,\xi\rangle .
\end{align*}
For a continuous $w$ the distributional right Djrbashian-Caputo derivative is
$\rCder{\beta}_t w:=-\partial_t I_{T-}^{1-\beta}(w-w(\cdot,T))$.
Thus $u(\cdot,T)=0$ yields
$\mathcal D^\beta\check u=\rCder{\beta}_t u$ in $\mathcal D'(Q_R)$.

 The translated-test argument in that proof gives, on each compact subset of $Q_R$ for all sufficiently large $n$,
\[
 \mathcal D^\beta U_n-\frac{\sigma^2}{2}(U_n)_{gg}\le C_R,
 \quad
 \mathcal D^\beta V_n-\frac{\sigma^2}{2}(V_n)_{gg}\ge-C_R
\]
in the viscosity sense. The regularizations are locally Lipschitz and respectively semiconvex and semiconcave; see \cite[Proposition~3.3]{Namba18} for these standard properties. Alexandrov's theorem \cite[Theorem~A.2]{CIL92} gives the inequalities almost everywhere. The singular spatial second-derivative measures are nonnegative for $U_n$ and nonpositive for $V_n$, so the respective inequalities also hold in distributions.

Let $\varepsilon_0:=\sup_n\varepsilon_n$. Since $\check u\ge0$, using
the source $(g,t)$ in both convolutions and completing the square in
the spatial penalty gives, for all $t\ge0$,
\[
 0\le V_n(g,t)\le\check u(g,t)\le U_n(g,t)
 \le C(1+|g|)+\tfrac12 C^2\varepsilon_0 .
\]
Thus both sequences are uniformly bounded on every
$[-R',R']\times[0,\infty)$. They converge locally uniformly to $\check u.$
For $W_n=U_n$ or $V_n$ and $0<\beta<1$, the strip bound and the
integrable tail of $(\mathcal D^\beta)^*\xi$ give, by dominated convergence,
\[
 \langle\mathcal D^\beta W_n,\xi\rangle
 \longrightarrow \langle\mathcal D^\beta\check u,\xi\rangle
 =\langle\rCder{\beta}_t u,\xi\rangle .
\]
The spatial derivatives pass to the limit by testing against $\xi_{gg}.$
For $\beta=1$ the time term passes to the limit by testing against
$\partial_t\xi$. Taking $\xi\ge0$ in the upper inequality for $U_n$
and the lower inequality for $V_n$ proves
\eqref{eq:time-heat-bounds} in $\mathcal D'(Q_R)$.

Let $\zeta\in C_c^\infty((-1,1))$ be a nonnegative smooth mollifier with $\int_{-1}^{1}\zeta(y)\,dy=1$. For $0<\varepsilon\le1$, define
\[
 \zeta_\varepsilon(y):=\varepsilon^{-1}\zeta(y/\varepsilon),
 \quad u_\varepsilon(g,t):=(\zeta_\varepsilon*_g u)(g,t)=\int_{-\varepsilon}^{\varepsilon}\zeta_\varepsilon(y)u(g-y,t)\,dy.
\]
Since $u$ is $L$-Lipschitz in $g$, $\|u_g\|_{L^\infty(\R\times(0,T))}\le L$. Standard convolution estimates and integration by parts give, for all $(g,t)\in\R\times[0,T]$ and $0<\varepsilon\le1$,
\begin{equation}\label{eq:moll-lip-bounds}
 |u_\varepsilon(g,t)-u(g,t)|\le C_1 L\varepsilon,\quad
 |(u_\varepsilon)_{gg}(g,t)|=|(\zeta_\varepsilon' *_g u_g)(g,t)|\le C_2 L\varepsilon^{-1},
\end{equation}
where $C_1:=\int_{-1}^{1}|y|\zeta(y)\,dy$ and $C_2:=\|\zeta'\|_{L^1(-1,1)}$.

Convolving \eqref{eq:time-heat-bounds} in space with $\zeta_\varepsilon\ge0$ shows that for $|g|\le R$, $u_\varepsilon(g,\cdot)$ satisfies $\rCder{\beta}_t u_\varepsilon(g,\cdot)=f_\varepsilon(g,\cdot)$ in $\mathcal D'(0,T)$, where by \eqref{eq:moll-lip-bounds} the function $f_\varepsilon:=\mathcal L_\beta u_\varepsilon+\frac{\sigma^2}{2}(u_\varepsilon)_{gg}$ satisfies 
$$
\|f_\varepsilon(g,\cdot)\|_{L^\infty(0,T)}\le C_R+\frac{\sigma^2}{2}C_2 L\varepsilon^{-1}\le \widetilde C_R\varepsilon^{-1}.
$$ 
Since $u_\varepsilon(g,\cdot)\in C([0,T])$ with $u_\varepsilon(g,T)=0$, its distributional fractional and Riemann--Liouville derivatives agree. For $0<\beta<1$, writing $I_{T-}^{1-\beta}$ for the right Riemann--Liouville fractional integral, we therefore have
\[
 -\partial_t\bigl(I_{T-}^{1-\beta}u_\varepsilon(g,\cdot)\bigr)
 =f_\varepsilon(g,\cdot)\in L^\infty(0,T).
\]
Therefore, $I_{T-}^{1-\beta}u_\varepsilon(g,\cdot)\in W^{1,\infty}(0,T)$, and its terminal trace is zero by boundedness of $u_\varepsilon$. The inversion formula \cite[Lemma~2.6(b)]{kilbas2006theory} thus yields, by continuity (and by ordinary integration when $\beta=1$),
\[
 u_\varepsilon(g,t)=\frac{1}{\Gamma(\beta)}\int_t^T(q-t)^{\beta-1}f_\varepsilon(g,q)\,dq,
 \quad t\in[0,T].
\]
By the standard $\beta$-H\"older estimate for Riemann--Liouville fractional integrals on $L^\infty(0,T)$ \cite[Section~3.1, Corollary~2 and the following estimate, p.~56]{SKM93}, applied after the time reflection $t\mapsto T-t$ (and by ordinary integration for $\beta=1$),
\[
 |u_\varepsilon(g,t)-u_\varepsilon(g,s)|
 \le\frac{2\|f_\varepsilon(g,\cdot)\|_{L^\infty(0,T)}}{\Gamma(1+\beta)}|t-s|^\beta
 \le C_R'\varepsilon^{-1}|t-s|^\beta.
\]
Combining this bound with the first estimate in \eqref{eq:moll-lip-bounds} via the triangle inequality gives, for $|g|\le R$ and $0<\varepsilon\le1$,
\begin{align*}
 |u(g,t)-u(g,s)|
 &\le |u(g,t)-u_\varepsilon(g,t)|+|u_\varepsilon(g,t)-u_\varepsilon(g,s)|+|u_\varepsilon(g,s)-u(g,s)|\\
 &\le C_R''\left(\varepsilon+\frac{|t-s|^\beta}{\varepsilon}\right).
\end{align*}
Choosing $\varepsilon=|t-s|^{\beta/2}$ for $0<|t-s|\le1$ (and using the boundedness of $u$ on $[-R,R]\times[0,T]$ for $|t-s|>1$) establishes \eqref{eq:temporal-holder}.
\end{proof}

\section{Connection to Optimal Switching Problem}\label{sec:switching}
Now, we consider the optimal switching problems that arise in the study of speculative financial bubbles. 
We  note that the spatial part of the operator is symmetric; that is, if we define the spatial reflection $\hat{u}(g,t) = u(-g,t)$, the operator satisfies
\begin{equation*}
F[\hat{u}](g,t) = \rCder{\beta}_{t}u(-g,t) + \rho g (-u_g(-g,t)) - \tfrac{\sigma^2}{2}u_{gg}(-g,t) + r u(-g,t) = F[u](-g,t).
\end{equation*}

\subsection{Related switching problem}

To illustrate the connection between the asset bubble and optimal switching problems, we  follow \cite[Section $5$]{barkhudaryan2022}. Let

\begin{equation*}
\psi_1(g,t) \equiv \psi(g,t), \quad \psi_2(g,t) \equiv \psi(-g,t).
\end{equation*}

The original problem can be rewritten as a weakly coupled switching system. By carefully mapping the signs to standard optimal switching theory, the system becomes

\begin{equation}\label{eq:switching}
\begin{cases}
\min\Bigl\{F[u_1](g,t),\;\; u_1 - u_2 + (-\psi_1)\Bigr\} = 0, \\
\min\Bigl\{F[u_2](g,t),\;\; u_2 - u_1 + (-\psi_2)\Bigr\} = 0,
\end{cases}
\quad (g,t)\in\R\times[0,T),
\end{equation}
with $u_1(g,t) = u(g,t),\;u_2(g,t) = u(-g,t)$, and the terminal conditions $u_1(g,T) = 0,\; u_2(g,T) = 0$. 

In the standard optimal switching literature, the terms added to the difference $u_i - u_j$ represent the transition costs $\Psi_{ij}$. Thus, our effective transition costs are $\Psi_{12} = -\psi_1$ and $\Psi_{21} = -\psi_2$.

For an optimal switching problem to be well-defined, it must satisfy a non-profit loop condition (i.e., $\Psi_{12} + \Psi_{21} \ge 0$), meaning agents cannot make arbitrary gains by continuously looping between states. In our specific case, the complement condition is given as
\begin{equation*}
  \psi(g,t)+\psi(-g,t)=-2c<0.
\end{equation*}
This implies that $\psi_1 + \psi_2 = -2c$. Translating this to our effective transition costs yields
\begin{equation*}
\Psi_{12} + \Psi_{21} = -\psi_1 - \psi_2 = -(-2c) = 2c > 0.
\end{equation*}

\begin{theorem}\label{thm:switching}
Because the effective switching costs satisfy the strict no-profit loop condition ($\Psi_{12} + \Psi_{21} = 2c > 0$), the fractional bubble problem \eqref{eq:pde}--\eqref{eq:operator} and the switching system \eqref{eq:switching} are equivalent in the following sense
\begin{enumerate}
    \item[(i)] If $u$ is the unique viscosity solution of \eqref{eq:pde}, then 
    $$
    (u_1(g,t), u_2(g,t)) := (u(g,t), u(-g,t))
    $$ 
    is a viscosity solution of the switching system \eqref{eq:switching}.
    \item[(ii)] Conversely, if $(u_1, u_2)$ is a viscosity solution of \eqref{eq:switching} in the linear-growth class, then by reflection symmetry of $F$ and $\psi_2(g,t)=\psi_1(-g,t)$, the pair $(u_2(-g,t), u_1(-g,t))$ also solves \eqref{eq:switching}. Uniqueness for \eqref{eq:switching} under the strict no-profit loop condition then enforces $u_2(g,t) = u_1(-g,t)$, and $u := u_1$ solves the fractional bubble problem \eqref{eq:pde}.
\end{enumerate}
\end{theorem}
\begin{proof}
    The proof can be done by using the arguments in the proof of Theorem $5.1$ in \cite{barkhudaryan2022}.
Let $\widehat w(g,t):=w(-g,t)$ and interpret \eqref{eq:switching}
componentwise as in Definition~\ref{def:visc}. Since $F$ commutes with
reflection and $\psi_2=\widehat{\psi_1}$, every solution $u$ of
\eqref{eq:pde} gives the switching pair $(u,\widehat u)$.
Conversely, for any switching solution of uniform linear growth, testing
the component attaining each extremum shows that
$(\max\{u_1,\widehat u_2\})^*$ is a scalar subsolution and
$(\min\{u_1,\widehat u_2\})_*$ a scalar supersolution. Here we recall that 
\(w^*\) stands for the upper semicontinuous envelope of \(w\), and \(w_*\) is its lower semicontinuous envelope.
The comparison principle (Theorem~\ref{thm:comp}) therefore forces $u_1=\widehat u_2$, and the
first equation of \eqref{eq:switching} reduces to \eqref{eq:pde}.
\end{proof}

\section*{Acknowledgments}

\noindent  

\medskip
\noindent {\bf Funding:}
A. Arakelyan was supported by the Higher Education and Science Committee of the MESCS RA under Research Project No. 24IRF-1A001.
R. Barkhudaryan  was supported by the Higher Education and Science Committee of the MESCS RA under Research Project No. 25RL-1A040.
L. El Hajj was supported by the Qatar Foundation through Carnegie Mellon University in Qatar’s Accelerator Research Funding program.  H. Shahgholian was supported by Swedish Research Council, grant No. 2025-03740.  
\medskip

\noindent{\bf AI Usage Statement:}
Artificial intelligence (AI) tools were used solely to assist with proofreading, refining language clarity, and verifying select mathematical computations during the preparation of this manuscript. All original intellectual contributions, core analyses, interpretation of results, and final writing are the sole responsibility of the authors.

\medskip

\noindent {\bf  Data availability statement:} All data needed are contained in the manuscript.

\medskip
\noindent {\bf Conflicts of interests/Competing interests:} The authors declare that there are no financial, competing or conflict of interests.

\bibliographystyle{acm} 
\bibliography{bubble-fBm}

\end{document}